\documentclass[matzei,envcountsame,draft]{svjour}

\usepackage[latin1]{inputenc}%
\usepackage{amsmath,amssymb,latexsym}%

	\def\makeheadbox{{%
	\hbox to0pt{\vbox{\baselineskip=10dd\hrule\hbox
	to\hsize{\vrule\kern3pt\vbox{\kern3pt
	\hbox{January 24th, 2006}
	\kern3pt}\hfil\kern3pt\vrule}\hrule}%
	\hss}}}

\def\Romannumeral#1 {\begingroup{\uppercase\expandafter{\romannumeral#1}}\endgroup}% 
\def\pf#1{{\def\temp{#1} 
              \ifx\temp\empty 
                  \noindent\slshape\textbf{Proof.\ }
              \else 
                  \noindent\slshape\textbf{Proof of\ #1.\ }
              \fi}}
\def\qed{\hspace*{\fill}$\square$}
\newcommand{\im}{\operatorname{im}}%
\newcommand{\rank}{\operatorname{rank}}%
\newcommand{\Hom}{\operatorname{Hom}}% 
\def\Stab#1_#2.{#2_{\{#1\}}}%
\def\Fix#1_#2.{{#2_{#1}}}%
\def\Aut#1_#2.{\mathsf{Aut}_{#2}(#1)}%
\def\Isom#1_#2.{\mathsf{Isom}_{#2}(#1)}%
\newcommand{\CAT}[1]{\mathsf{CAT}(#1)}%
\def\bs#1.{
              \def\temp{#1} 
              \ifx\temp\empty 
                   \mathcal{B}
              \else
                   \mathcal{B}(#1)
              \fi
}

\renewcommand{\emph}{\textbf}%

\newcommand{\frk}{\mathop{\text{\rm flat-rk}}}
\newcommand{\ark}{\mathop\text{\rm alg-rk}}
\newcommand{\rk}{\mathop\mathrm{rk}}

\spnewtheorem{algorithm}[theorem]{Algorithm}{\bf}{\sf}

\begin{document}

\title{Flat rank of automorphism groups of buildings}

\author{
Udo Baumgartner\inst{1}, Bertrand R\'emy\inst{2}  \and George A. Willis\inst{1}
}                     

\institute{
School of Mathematical and Physical Sciences, 
The University of Newcastle, University Drive, Building V, 
Callaghan, NSW 2308, Australia\\
\email {Udo.Baumgartner@newcastle.edu.au, George.Willis@newcastle.edu.au}
\and
 Institut Camille Jordan, UMR 5208 CNRS / Lyon 1, Universit\'e Claude Bernard Lyon 1, 
21 avenue Claude Bernard, 69622 Villeurbanne cedex, France\\ 
\email{remy@math.univ-lyon1.fr}
}

\date{}

\maketitle

\begin{abstract}

The flat rank of a totally disconnected locally compact group $G$, denoted $\frk(G)$, is an invariant of the topological group structure of $G$. 
It is defined thanks to a natural distance on the space of compact open subgroups of $G$. 
For a topological Kac-Moody group $G$ with Weyl group $W$, we derive the inequalities: 
$\ark(W)\le \frk(G)\le \rk(|W|_0)$. 
Here, $\ark(W)$ is the maximal $\mathbb{Z}$-rank of abelian subgroups of $W$, and $\rk(|W|_0)$ is the maximal dimension of isometrically embedded flats in the $\CAT0$-realization $|W|_0$. 
% I removed Davis, because I would say Davis for buildings (and Moussong for Coxeter groups)
We can prove these inequalities under weaker assumptions. 
We also show that for any integer $n \geq 1$ there is a topologically simple, compactly generated, locally compact, totally disconnected group $G$, with $\frk(G)=n$ and which is not linear. 
% which group would you take for n=0? 

\keywords {
totally disconnected group, flat rank, automorphism group, scale function, twin building, strong transitivity, Kac-Moody group.
}

\subclass{
22D05, 22D45 (primary). 
%	20E36, % (General theorems concerning automorphisms of groups), 
%	20E08, % (groups acting on trees) 
%	20G25  %(Linear algebraic groups over local fields and their integers)
%	(secondary)
}

\end{abstract}

\section*{Introduction}\label{sec:intro}
The general structure theory of locally compact groups is a well-established topic in mathematics. 
One of its main achievements is the solution to Hilbert's 5th problem on characterizing Lie groups. 
The general structure results \cite{lcG.trafo} are still used in recent works. 
For instance, the No Small Subgroup Theorem  \cite[4.2]{lcG.trafo}
is used in 
Gromov's characterization 
of finitely generated groups 
of polynomial growth \cite{G(poly-growth)+expand-maps}. 
More recently, 
the theory was used 
in rigidity problems for discrete groups, 
in order to 
attach suitable boundaries 
to quite general topological groups \cite{cont-bound-cohom+app(rig-Th)}. 

Simultaneously to these applications, 
recent years have seen 
substantial progress 
in extending known results 
about connected, locally compact groups 
to arbitrary locally compact groups. 
For example 
in the theory of random walks on groups 
it was shown in \cite{lcG.conc-funct} 
that 
the concentration functions 
for an irreducible probability measure 
on a non-compact group 
converge to $0$, 
while \cite{lcG:dense-orbits(Zd-act(aut))} 
contains the classification of 
ergodic $\mathbb{Z}^d$-actions 
on a locally compact group 
by automorphisms.  
%	The paper 
%	\cite{contrG+scales(AUT(tdlcG))}
%	finally 
%	describes contraction groups 
%direction(aut(tdlcG))%
%

This progress 
has been due to 
structure theorems 
for totally disconnected, locally compact groups 
established in 
\cite{tdlcG.structure}, and further advanced in \cite{furtherP(s(tdG))} and \cite{tidy<:commAut(tdlcG)}. 

The study of 
particular classes of examples 
has played an important role 
in informing  
further developments 
of the structure theory 
of totally disconnected, locally compact groups,  
beginning with 
the study of 
the classes of $p$-adic Lie groups 
\cite{scaleF(p-adic.Lie),uniscalar.padicLieG} 
and 
automorphism groups of graphs \cite{struc(tdlcG-graphs+permutations)}. 
This paper 
starts  
the examination of the topological invariants 
of totally disconnected groups 
which are 
closed automorphism groups of buildings with sufficiently transitive actions. 

%For reasons to be explained shortly,  
Topological Kac-Moody groups 
form a subclass 
of the latter class of groups, 
which  
is of particular interest to us.  
% 
%	Topological Kac-Moody groups 
%	form a 
%	large,  
%	relatively well behaved 
%	subclass of 
%	one of the natural classes 
%	of totally disconnected, locally compact groups,  
%	the closed subgroups of 
%	automorphisms groups of 
%	locally finite cell complexes, 
%	and they are 
%	sufficiently different from 
%	the classical case of $p$-adic Lie groups 
%	to expect new insight.  
%
From a combinatorial viewpoint 
\cite[1.C]{tG(Kac-Moody-type)+rangle-twinnings+:lattices} 
topological Kac-Moody groups 
generalize semi-simple algebraic groups 
and therefore 
should be expected to inherit 
some of the properties of linear groups. 
For instance, their Tits system structure and the virtual pro-$p$-ness of their maximal compact subgroups are used to prove their topological simplicity 
\cite[2.A.1]{top-simpl&com-superrig&non-lin(Kac-MoodyG)}; 
these properties are well known in the algebraic case.   
On the other hand, it is known that some of these groups are non-linear 
\cite[4.C.1]{top-simpl&com-superrig&non-lin(Kac-MoodyG)}. 
This may imply that 
the topological invariants 
of topological Kac-Moody groups 
differ substantially from 
those of 
algebraic groups over local fields 
in some important aspects. 
In this context 
we mention 
a challenging question, 
which also 
motivates our interest in 
topological Kac-Moody groups,
namely 
whether 
a classification of 
topologically simple, 
compactly generated, 
totally disconnected, 
locally compact groups 
is a reasonable goal.  

In this paper 
we focus on 
the most basic topological invariant 
of topological Kac-Moody groups $G$, 
the \emph{flat rank} of $G$, denoted $\frk(G)$. 
%
% Roughly speaking, 
This rank is defined using the space $\bs G.$ of compact open subgroups of $G$. 
This space is endowed with a natural distance: for $V,W\!\in\!\bs G.$, the numbers
$d(V,W)=\log\bigl(|V:V\cap W|\cdot |W:V\cap W|\bigr)$ define a discrete metric, for which conjugations in $G$ are isometries. 
A subgroup $O\leqslant G$ is called \emph{tidy} for an element $g\!\in\!G$ if it minimizes the displacement function of $g$, and a subgroup $H$ is called \emph{flat} if all its elements have a common tidy subgroup. 
A flat subgroup $H$ has a natural abelian quotient, whose rank is called its flat rank. 
Finally, the flat rank of $G$ is the supremum of the flat ranks of the flat subgroups of $G$. 
For details, we refer to Subsection \ref{subsec:fundamentals-tidy}. 
The two main results of the paper provide an upper and a lower bound for the flat rank of a sufficiently transitive automorphism group $G$ of a locally finite building. 

A summary of the main results about the upper bound for $\frk(G)$ is given by the following statement. 
The assumptions made 
in this theorem 
are satisfied by 
topological Kac-Moody groups
% ; see 
(Subsection~\ref{subsec:topAutG(buildings)}).

\vspace{3mm}
\noindent\emph{Theorem A.}~
{\em 
Let $(\mathcal{C},S)$ be a locally finite building with Weyl group $W$. 
%We 
Denote by 
$\delta$ and by $X$ 
the $W$-distance function and the $\CAT0$-realization of $(\mathcal{C},S)$, respectively. 
Let $G$ be a closed subgroup of the group of automorphisms of $(\mathcal{C},S)$. 
%We 
Assume that the $G$-action is transitive on ordered pairs of chambers at given $\delta$-distance. 
Then the following statements hold. 
\begin{itemize}
\item[{\rm (\romannumeral1)}]
The map $\varphi\colon X \to \bs G.$ mapping a point to its stabilizer, is a quasi-isometric embedding.
\item[{\rm (\romannumeral2)}]
For any point $x\!\in\!X$, the image of the orbit map $g \mapsto g.x$, restricted to a flat subgroup of flat rank $n$ in $G$, is an $n$-dimensional quasi-flat of $X$. 
\item[{\rm (\romannumeral3)}]
We have: $\frk(G)\le \rk(X)$. 
\item[{\rm (\romannumeral4)}]
If $X$ contains an $n$-dimensional flat, so does any of its apartments. 
\end{itemize}
\noindent As a consequence, we obtain: $\frk(G)\le \rk(|W|_0)$, where $\rk(|W|_0)$ is the maximal dimension of flats in the $\CAT0$-realization $|W|_0$. 
}
\vspace{3mm}

The strategy for the proof of Theorem~A 
is as follows.   
The inequality $\frk(G)\le \rk(|W|_0)$ 
is obtained from 
the statements (\romannumeral1)--(\romannumeral4),  
which are proven 
in the  listed order under weaker hypotheses. 
Statement 
(i) is part of Theorem \ref{thm:Main-delta-2-transitive} --- 
usually called the Comparison Theorem in this paper, 
(ii) is in fact proved in Proposition \ref{prop:nice-action=>flat-subgr->quasi-flat} under (i) as assumption, (iii) is a formal consequence of results of Kleiner's, and 
(iv) is Proposition \ref{prop:rk(building)=rk(apartment)}. 

The main result 
about the lower bound for $\frk(G)$ 
is the second half of Theorem \ref{thm:lifting-free_ab<(W)}, 
which we reproduce here as Theorem~B. 
The class of groups 
satisfying the assumptions of Theorem~B 
is contained in 
the class of groups 
satisfying the assumptions of Theorem~A 
and contains 
all topological Kac-Moody groups 
(Subsection~\ref{subsec:topAutG(buildings)}). 
% (see again subsection~\ref{subsec:topAutG(buildings)}). 

\vspace{3mm}
\noindent\emph{Theorem B.}~
{\em 
Let $G$ be a group with a locally finite twin root datum of associated Weyl group $W$. 
We denote by $\overline{G}$ the geometric completion of $G$, i.e. the closure of the $G$-action in the full automorphism group of the positive building of $G$. 
Let $A$ be an abelian subgroup of $W$. 
Then $A$ lifts to a flat subgroup $\widetilde{A}$ of $\overline{G}$ such that $\frk(\widetilde{A})=\rank_\mathbb{Z}(A)$.\\ 
As a consequence, we obtain: $\ark(W)\le \frk(G)$, where $\ark(W)$ is the maximal $\mathbb{Z}$-rank of abelian subgroups of $W$. 
}
\vspace{3mm}

The first half of Theorem \ref{thm:lifting-free_ab<(W)} (not stated here) asserts that the flat rank of the rational points of a semi-simple group over a local field $k$ coincides with the algebraic $k$-rank of this group. 
This is enough to exhibit topologically simple, compactly generated, locally compact, totally disconnected groups of arbitrary flat rank $d \geq 1$, e.g. by taking the sequence 
$\bigl({\rm PSL}_{d+1}(\mathbb{Q}_p)\bigr)_{d\geq 1}$. 
Theorem~C (Theorem \ref{thm:E(non-linear-simple-tdlcG(any_rank))} in the text)
% will enable 
enables us to 
exhibit a sequence of {\it non-linear}~groups with the same properties:  

\vspace{3mm}
\noindent\emph{Theorem C.}~
{\em 
For every integer $n \ge 1$ there is a non-linear, topologically simple, compactly generated, 
locally compact, totally disconnected group of flat rank  $n$. 
}
\vspace{3mm}

These examples are provided by Kac-Moody groups. 
The combinatorial data leading to these groups are obtained by gluing a hyperbolic Coxeter diagram arising from a non-linear Kac-Moody group, together with an affine diagram which ensures the existence of a sufficiently large abelian group in the resulting Weyl group. 

Let us finish this introduction with a conjecture. 
In Theorems A and B,
the upper and lower bound on $\frk(G)$ 
depend only on the associated Weyl group $W$.
We think that the bounds are equal;  
if this is indeed so, 
the present paper computes 
the exact rank 
of any geometric completion 
of a finitely generated Kac-Moody group 
thanks to 
a result in Daan Krammer's PhD thesis, 
\cite[Theorem~6.8.3]{thesisKrammer}; 
see Subsection~\ref{subsec:quotingKrammer}. 
%provides a combinatorial control on the abelian subgroups of a Coxeter group. 
%However it is not clear to us that we can use it to obtain information on the maximal flats in $|W|_0$. 

\vspace{3mm}
\noindent\emph{Conjecture.}~
{\em 
Let $W$ be a finitely generated Coxeter group.
Then we have: $\rk(|W|_0)=\ark(W)$, where $\rk(|W|_0)$ is the maximal dimension of flats in the $\CAT0$-realization of the Coxeter complex of $W$ and $\ark(W)$ is the maximal rank of free abelian subgroups of $W$.
}
\vspace{3mm}

This conjecture was checked by Fr\'ed\'eric Haglund 
in the case of right-angled Coxeter groups \cite{Haglund}. 
By Bieberbach theorem, 
to prove this conjecture 
it is enough to show 
the existence of  
a ``periodic'' flat 
of dimension equal to $\rk(|W|_0)$ 
in $|W|_0$; 
i. e. one  
which admits a cocompact action 
by a subgroup of $W$.
% 

%Even though 
%it is of no importance 
%for the topics discussed in this paper 
The following, 
more general,  
conjecture, 
seems to be the natural framework 
for questions of this kind: 
\textit{If 
a group $G$ acts cocompactly on a proper $\CAT0$-space $X$,  
then 
the rank of $X$ equals 
the maximal rank of an abelian subgroup of $G$.
}
Some singular cases 
\cite{BalBrinGD,BalBrinIHES} 
of this generalization  
as well as 
the smooth analytic case
\cite{BBSsmooth,BanSchro} 
have been proved.
\vspace{3mm} 

\textsc{Organization of the paper.}~
Section \ref{sec:fundamentals} is devoted to recalling basic facts on buildings (chamber systems), on the combinatorial approach to Kac-Moody groups (twin root data), and on the structure theory of totally disconnected locally compact groups (flat rank). 
We prove the upper bound inequality in Section~\ref{sec:upper_bound}, and the lower bound inequality in Section~\ref{sec:lower_bound}.
In Section~\ref{sec:cases_equality}, we exhibit a family of non-linear groups of any desired positive flat rank; we also discuss applications of the flat rank to isomorphism problems. 

\vspace{3mm} 

\textsc{Acknowledgements:}~
We thank Pierre--Emmanuel Caprace for making us aware of Daan Krammer's work, as well as Werner Ballmann, Michael Davis and Fr\'ed\'eric Haglund for useful conversations about the above conjecture.

% SECTION 1

\section{Buildings and totally disconnected groups}\label{sec:fundamentals}

\subsection{Buildings and their automorphism groups}\label{subsec:fundamentals-buildings}

\subsubsection{Buildings as chamber systems.}
In this paper 
a building is a chamber system, 
throughout denoted $\mathcal{C}$, 
together with a distance function 
with values in a Coxeter group.  
A chamber system 
is a set $\mathcal{C}$, 
called \emph{the set of chambers},   
together with a family $S$ of partitions of $\mathcal{C}$, 
called \emph{the adjacency relations}. 
Each element $s$ of $S$ 
defines an equivalence relation 
which we will not distinguish from $s$. 
For 
each 
$s$ in $S$ 
the equivalence classes 
of $s$ 
should be thought of 
as the set of chambers 
sharing a fixed `face' 
of `color' $s$.

A finite sequence of chambers 
such that consecutive members are adjacent 
(that is, contained in some adjacency relation) 
is called a \emph{gallery}. 
A gallery is said to \emph{join} its first and last term. 
A gallery is called \emph{non-stuttering} 
if consecutive members are different. 
For every non-stuttering gallery $(c_0,\ldots,c_n)$, 
any word $s_1\cdots s_n$ in the free monoid $S^*$ on $S$  
such that $c_{j-1}$ and $c_j$ are $s_j$-adjacent 
for all $1\le j\le n$, is called 
a \emph{type} of the gallery $(c_0,\ldots,c_n)$ 
(a gallery 
does not necessarily have 
a unique type). 
A gallery having a type 
contained in 
the submonoid of $S^*$ 
generated by a subset $T$ of $S$ 
is called a \emph{$T$-gallery}. 
A maximal subset of chambers 
which can be joined by a $T$-gallery 
is called a \emph{$T$-residue}.  
We call a chamber system $(\mathcal{C},S)$ \emph{locally finite} 
if $S$ and all $\{s\}$-residues for $s\!\in\!S$ are finite. 

A permutation of the underlying set $\mathcal{C}$ 
of a chamber system $(\mathcal{C},S)$
is said to be an \emph{automorphism} of $(\mathcal{C},S)$ 
if it induces a permutation of $S$. 
The group of all automorphisms of $(\mathcal{C},S)$ 
will be denoted $\Aut{\mathcal{C}, S}_.$. 
An automorphism of  $(\mathcal{C},S)$
is said to be \emph{type-preserving} 
if it induces the identity permutation of $S$. 
The group of all type-preserving automorphisms of $(\mathcal{C},S)$ 
will be denoted $\Aut{\mathcal{C}, S}_0.$. 

Each Coxeter system $(W,S)$ gives rise to a chamber system: 
its set of chambers is $W$, 
and for $s\!\in\!S$, we say that 
$w$ and $w'$ are $s$-adjacent 
if and only if 
$w'\in \{w,ws\}$.  
A word $f$ in the free monoid on $S$ 
is called \emph{reduced} 
if it has minimal length among all such words representing 
their product $s_f$ 
as an element of $W$. 
If $(W,S)$ is a Coxeter system 
and $T$ is a subset of $S$ 
then the subgroup of $W$ generated by $T$ 
is called a \emph{special subgroup} 
and is denoted $W_T$. 
A subset $T$ of $S$ 
is called \emph{spherical} 
if $W_T$ is finite. 

Let $(W,S)$ be a Coxeter system. 
A \emph{building} of type $(W,S)$ 
is a chamber system $\mathcal{C}$ 
with adjacency relations indexed by the elements of $S$,  
each consisting of sets containing at least $2$ elements. 
We also require the existence of a \emph{$W$-distance function} 
$\delta\colon \mathcal{C}\times\mathcal{C} \to W$ 
such that 
whenever $f$ is a reduced word on $S$,  
then for chambers $x$ and $y$ 
we have $\delta(x,y)= s_f$ 
if and only if 
there is a (non-stuttering) 
gallery of type $f$ joining $x$ to $y$. 
In a building 
a non-stuttering gallery 
has a unique type. 
Other 
basic properties of $W$-distance functions 
can be found in \cite[3.1]{PiM7}.  

\subsubsection{Non-positively curved realization of a building.}
A chamber system can be realized 
as a topological space so that each chamber is homeomorphic to 
a model space $X$, and adjacency 
of chambers is represented by 
them sharing a preassigned subspace of $X$ 
as a common `face'. 
We now explain 
a very flexible way to do this 
for chamber systems 
which are buildings. 
Moussong attributes the method to Vinberg \cite {discrete-linG_gen(reflections)}. 
We follow Davis' exposition \cite{buildings=CAT0}. 

Let $(\mathcal{C},S)$ be a building of type $(W,S)$. 
We start out with 
a topological space $X$, 
which will be our model for a chamber, 
and a family of closed subspaces $(X_s)_{s \in S}$, 
which will be our supply of `faces'. 
The pair $(X,(X_s)_{s \in S})$ will be called a \emph{model space}. 
For each point 
$x$ in $X$ 
we define a subset $S(x)$ of $S$ by setting 
$S(x):=\{s\!\in\!S\colon x\in X_s\}$. 
Further, 
we define an equivalence relation $\sim$ 
on the set $\mathcal{C}\times X$ by 
$(c,x)\sim (c',x')$ if and only if $x=x'$ and $\delta(c,c')\in W_{S(x)}$. 
The \textbf{$X$-realization of $(\mathcal{C},S)$}, 
written $X(\mathcal{C})$, is 
the quotient space $(\mathcal{C}\times X)/\!\!\sim$, 
where $\mathcal{C}$ carries the discrete topology. 

%	It can be desirable 
%	to choose a model space $(X,(X_s)_{s \in S})$ 
%	such that 
%	the images of two chambers 
%	$c$ and $c'$ 
%	under the map 
%	$\mathcal{C}\to \mathcal{C}\times X\to X(\mathcal{C})$ 
%	do not share any points 
%	even though 
%	$c$ and $c'$ 
%	lie in the same residue. 
%	Indeed, 
%	an important feature of 
%	the Davis realization of a building,  
%	which will arise from a specific choice of model space 
%	to be explained below, 
%	is that 
%	the geometric realizations 
%	of two chambers 
%	have common points 
%	if and only if 
%	they lie in 
%	a common spherical residue. 
 
 If $(\mathcal{C},S)$ is a building, 
 and 
 $(X,(X_s)_{s \in S})$ is a model space,  
 then 
 any type-preserving automorphism of $(\mathcal{C},S)$ 
 induces a homeomorphism of $X(\mathcal{C})$ 
 via the induced permutation action on $\mathcal{C}\times X$. 
 Further, 
 this assignment 
 defines a homomorphism 
 from 
$\Aut{\mathcal{C}, S}_0.$ 
into the group of homeomorphisms of $X(\mathcal{C})$. 
This homomorphism is injective if  $X\setminus \bigcup_{s \in S} X_s \neq\varnothing$.  
Automorphisms of $(\mathcal{C},S)$ 
which are not type-preserving 
will not induce homeomorphisms of $X(\mathcal{C})$ 
unless the model space $(X,(X_s)_{s \in S})$ 
admits symmetries 
realizing the possible type permutations. 
We will spell out 
appropriate conditions below 
for specific choices of model spaces. 

We now introduce the model spaces $(X,(X_s)_{s \in S})$ 
which define  
the Davis realization 
of a building $(\mathcal{C},S)$ 
leading to a $\CAT0$-structure on $X(\mathcal{C})$. 
A variant of it, available for a subclass of buildings, 
is used by Moussong to define a $\CAT{-1}$-structure on $X(\mathcal{C})$ \cite{hyp-CoxG}. 
{\it We assume from now on that 
$(\mathcal{C},S)$ is a building of finite rank, i.e. $S$ is finite}. 
The model spaces $X$ in both cases 
are metric simplicial complexes 
(with a family of subcomplexes $(X_s)_{s \in S}$) 
with the same underlying abstract simplicial complex,  
namely 
the flag complex of the poset of spherical subsets of $S$ 
ordered by inclusion. 
For $s\!\in\!S$ 
the subcomplex $X_s$ is 
the union of all 
chains starting with the set $\{s\}$ 
(and all their sub-chains). 
This model space 
always supports a natural piecewise Euclidean structure \cite{buildings=CAT0} 
as well as a piecewise hyperbolic structure \cite[section~13]{hyp-CoxG} 

Our assumption that $S$ is a finite set 
implies that 
$X$ is a finite complex; 
in particular $X(\mathcal{C})$ has only finitely many cells up to isometry. 
Therefore Bridson's theorem \cite[I.7.50]{GMW319}
implies that $X(\mathcal{C})$ with the path metric is a complete geodesic space. 
Moreover $X$ has finite diameter since only compact simplices are used for the hyperbolic structure.  
Suppose that all $s$-residues of $(\mathcal{C},S)$ are finite for $s\!\in\!S$. 
Then, because of the way we defined the family of subspaces $(X_s)_{s \in S}$ encoding the faces, 
$X(\mathcal{C})$ is locally finite.   
The geometric realization of 
a Coxeter complex based on model spaces for the Davis realization is $\CAT0$. 
The Moussong realization of a Coxeter complex is $\CAT{-1}$ if and only if the Coxeter group is Gromov-hyperbolic. 
Both results are derived in \cite{hyp-CoxG}. 
Using retractions onto apartments one shows that analogous results hold for buildings whose Coxeter group is of the appropriate type \cite[section~11]{buildings=CAT0}. 

% I am not sure it is a good idea to call Moussong what is CAT(-1) 
% and Davis what is simply CAT(0)
% from my viewpoint, the difference lies between thin and thick buildings

For both the Davis and Moussong realizations 
the map which assigns to a type-preserving automorphism 
of the building $(\mathcal{C},S)$ 
the self-map of $X(\mathcal{C})$ 
induced by the permutation of $\mathcal{C}$ 
defines a homomorphism 
from  $\Aut{\mathcal{C}, S}_0.$ 
into the group of simplicial isometries of $X(\mathcal{C})$, 
which we will denote by $\Isom X(\mathcal{C})_.$.   
The metric structures on the corresponding model spaces 
are in addition 
invariant under all diagram automorphisms of the 
Coxeter diagram of the building. 
Hence 
automorphisms of $(\mathcal{C},S)$ 
also induce simplicial isometries of $X(\mathcal{C})$ 
in both cases 
and the analogous result holds for 
automorphisms of $(\mathcal{C},S)$ 
which are not necessarily type-preserving. 
Since 
the vertex $\varnothing$ of $X$ 
is not contained in any of the subcomplexes $X_s$ for $s\!\in\!S$,  
these homomorphisms are injective. 
We denote the Davis and Moussong realizations of a building $(\mathcal{C},S)$ by $|\mathcal{C}|_0$ and $|\mathcal{C}|_{-1}$ respectively. 

The natural topology on the group of automorphisms 
of a chamber system is the permutation topology: 

\begin{definition}\label{def:permutation topology}
Suppose that 
a group $G$ 
acts on a set $M$. 
Denote 
the stabilizer 
of a subset $F$ of $M$ 
by $\Stab F_G.$. 
The \emph{permutation topology} 
on $G$ 
is the topology 
with 
the family $\{\Stab F_G.\colon  \text{$F$ a finite subset of $M$}\}$ 
as neighborhood base of the identity 
in $G$.  
\end{definition}

As we already noted,  
automorphisms of a chamber system 
can be viewed as 
permutations of the set of chambers. 
%Proposition~\ref{prop:Aut(chamber-sys)=Isom(geom-realization)} 
%allows us to identify the group of automorphisms 
%of a building with the group of isometries of its geometric realization. 
The permutation topology 
maps to another natural topology 
under the injection $\Aut{\mathcal{C}, S}_.\to \Isom{|\mathcal{C}|_\epsilon}_.$ 
for $\epsilon\in\{0,-1\}$.

\begin{lemma}\label{lem:perm-top->cp-open-top}
Let $(\mathcal{C},S)$ be a building 
with $S$ and all $s$-residues finite for $s\!\in\!S$. 
Then 
the permutation topology on $\Aut{\mathcal{C}, S}_.$ 
maps to the compact-open topology 
under the map $\Aut{\mathcal{C}, S}_.\to \Isom{|\mathcal{C}|_\epsilon}_.$ 
for $\epsilon\in\{0,-1\}$. 
\qed
\end{lemma}

\subsection{Topological automorphism groups of buildings}\label{subsec:topAutG(buildings)}
The examples of topological automorphism groups 
of buildings 
we are most interested in 
are Kac-Moody groups over finite fields. 
We will not define them 
and refer the reader to 
\cite[Subsection~3.6]{unique+pres(Kac-MoodyG(F))} 
and 
\cite[Section~9]{GKac-Moody-depl+presque-depl} 
for details instead. 
A Kac-Moody group 
over a finite field is an example of 
a group $G$ with twin root datum  $\bigl((U_a)_{a\in \Phi},H\bigr)$ 
(of type $(W,S)$; 
compare \cite[1.5.1]{GKac-Moody-depl+presque-depl} for the definition) 
such that 
all the root groups are finite.  
We will call a group 
which admits 
a twin root datum 
consisting of finite groups 
\emph{a group with a locally finite twin root datum}. 

Any group $G$ with locally finite twin root datum 
of type $(W,S)$ 
admits an action on a twin building 
(compare \cite[2.5.1]{GKac-Moody-depl+presque-depl} for the definition) 
of type $(W,S)$ 
having finite $s$-residues for all $s\!\in\!S$. 
Let $(\mathcal{C},S)$ be the positive twin 
(which is isomorphic to the negative twin). 
Its geometric realizations 
$|\mathcal{C}|_0$ and $|\mathcal{C}|_{-1}$ 
are locally finite. 
The action of $G$ on $(\mathcal{C},S)$ 
is strongly transitive in the sense of \cite[p. 56]{PiM7} 
(called strongly transitive with respect to an apartment in \cite[2.6.1]{GKac-Moody-depl+presque-depl}). 
The group $H$ 
is the fixator of an apartment $\mathcal{A}$ of $(\mathcal{C},S)$
with respect to this action of $G$. 
Hence 
we have a short exact sequence 
$1\to H\to N\to W\to 1$, 
where $N$ is the stabilizer of $\mathcal{A}$ 
(note that both BN-pairs are saturated). 

If $G$ is a group with a locally finite twin root datum $\bigl((U_a)_{a\in \Phi},H\bigr)$, then the $G$-actions on each building have a common kernel $K$. 
Moreover the root groups embed in $G/K$ and $\bigl((U_a)_{a \in \Phi}, H/K \bigr)$ is a twin root datum for $G/K$ with the same associated Coxeter system and twin buildings \cite[Lemma 1]{tG(Kac-Moody-type)+rangle-twinnings+:lattices}.
When $G$ is a Kac-Moody group over a finite field, we have $K=Z(G)$. 
The \emph{topological group associated to $G$}, 
denoted $\overline{G}$,  
is the closure of $G/K$ 
with respect to the topology 
on $\Aut{\mathcal{C}, S}_.$ 
defined in the previous subsection.

\subsection{Structure of totally disconnected, locally compact groups}\label{subsec:fundamentals-tidy}
The structure theory of totally disconnected, locally compact groups 
is based on 
the notions of tidy subgroup 
for an automorphism 
and the scale function. 
These notions 
were defined in \cite{tdlcG.structure} 
in terms of 
the topological dynamics of automorphisms 
and the definitions were reformulated in \cite{furtherP(s(tdG))}. 
We take the geometric approach 
to the theory 
as outlined in \cite{direction(aut(tdlcG))} 
and further elaborated on in \cite{G-flatrk_le1}. 

Let $G$ be a totally disconnected, locally compact group and let $\Aut G_.$ be the group of bicontinuous  automorphisms of $G$. 
We want to analyze the action of  subgroups of $\Aut G_.$ on $G$; 
we will be primarily interested in groups of inner automorphisms of $G$. 
To that end we consider the induced action of $\Aut G_.$ on the set  
\[
\bs G.:=\{V\colon \text{$V$ is a compact, open subgroup of $G$}\}\,.
\]

The function 
\[
d(V,W):=\log(|V\colon V\cap W|\cdot |W\colon W\cap V|)
\]
defines a metric on $\bs G.$ and $\Aut G_.$ acts by isometries on the discrete metric space $(\bs G.,d)$.
Let $\alpha$ be an automorphism of $G$. 
An element $O$ of $\bs G.$ 
is called \emph{tidy for $\alpha$} 
if the \emph{displacement function} of $\alpha$, 
denoted by $d_\alpha\colon \bs G.\to \mathbb{R}$ 
and defined by  $d_\alpha(V)=d(\alpha(V),V)$,  
attains its minimum at $O$. 
Since 
the set of values 
of the metric $d$ 
on $\bs G.$
is a well-ordered discrete subset of $\mathbb{R}$, 
every $\alpha\in\Aut G_.$ has a subgroup tidy for $\alpha$. 
Suppose that $O$ is tidy for $\alpha$. 
The integer 
\[
s_G(\alpha):= 
|\alpha(O)\colon \alpha(O)\cap O|\,, 
\]
which is also equal to 
$\min\bigl\{ |\alpha(V)\colon \alpha(V)\cap V|\colon V\in\bs G.\bigr\}$, 
is called \emph{the scale of $\alpha$}. 
An element of $\bs G.$ is called 
\emph{tidy for a subset} $\mathcal{M}$ of $\Aut G_.$ 
if and only if 
it is tidy for every element of $\mathcal{M}$. 
A subgroup $\mathcal{H}$ of $\Aut G_.$ 
is called \emph{flat} 
if and only if 
there is an element of $\bs G.$ 
which is tidy for $\mathcal{H}$. 
We will call a subgroup $H$ of $G$ flat 
if and only if 
the group of inner automorphisms 
induced by $H$ 
is flat. 

Later 
we will uncover implications 
of the flatness condition 
for groups acting 
in a nice way 
on $\CAT0$-spaces. 
They are based on 
the following properties of flat groups 
which hold 
for automorphism groups 
of general totally disconnected, locally compact groups.  
Suppose that 
$\mathcal{H}$ is a flat group of automorphisms.  
The set 
\[
\mathcal{H}(1):=\{\alpha\in\mathcal{H}\colon s_G(\alpha)=1=s_G(\alpha^{-1})\}
\]
is a normal subgroup of $\mathcal{H}$ 
and $\mathcal{H}/\mathcal{H}(1)$ is free abelian. 
The \emph{flat rank } of $\mathcal{H}$, 
denoted $\frk(\mathcal{H})$,  
is the $\mathbb{Z}$-rank of $\mathcal{H}/\mathcal{H}(1)$. 
If $\mathcal{A}$ is a group of automorphisms of 
the totally disconnected, locally compact group $G$ 
then its flat rank  
is defined to be 
the supremum of the flat ranks of all flat subgroups of $\mathcal{A}$. 
The flat rank  of 
the group $G$ itself 
is the flat rank  of the group of inner automorphisms of $G$.

If $\mathcal{H}$ is a flat group of automorphisms with $O$ tidy for $\mathcal{H}$ then by setting 
$\|\alpha \mathcal{H}(1)\|_\mathcal{H}:=d(\alpha(O),O)$, 
one defines a norm on $\mathcal{H}/\mathcal{H}(1)$ 
\cite[Lemma~11]{G-flatrk_le1}. 
That is, 
$\|\cdot\|_\mathcal{H}$ satisfies the axioms of a norm on a vector space 
with the exception 
that we restrict scalar multiplication to integers. 
This norm can be given explicitly in terms of a set of epimorphisms 
$\Phi(\mathcal{H},G)\subseteq \Hom(\mathcal{H},\mathbb{Z})$ 
of \emph{root functions} and a set of \emph{scaling factors} $s_\rho$; $\rho\in \Phi(\mathcal{H},G)$
associated to $\mathcal{H}$. 
In terms of these invariants of $\mathcal{H}$ the norm may be expressed as 
$
\|\alpha\mathcal{H}(1)\|_\mathcal{H}= \sum_{\rho\in\Phi} \log(s_\rho)\,|\rho(\alpha)|
$ 
\cite[Remark~13]{G-flatrk_le1}.
In particular the function $\|\cdot\|_\mathcal{H}$ extends to a norm in the ordinary sense on the vector space $\mathbb{R}\otimes \mathcal{H}/\mathcal{H}(1)$. 
For further information on flat groups see \cite{tidy<:commAut(tdlcG)}.

% SECTION 2

\section{Geometric rank as an upper bound}\label{sec:upper_bound}
In this section, we show that 
the flat rank  of 
a locally compact, strongly transitive group $G$ 
of automorphisms of a locally finite building, 
is bounded above by 
the geometric rank of the Weyl group. 
This inequality actually holds 
under the weaker assumption that 
the $G$-action is $\delta$-2-transitive: 

\begin{definition}\label{def:delta-2-transitive} 
Suppose the group $G$ 
acts on 
a building $(\mathcal{C},S)$ 
with Weyl group $W$ 
and $W$-valued distance function $\delta$. 
Then the action of $G$ 
on  $(\mathcal{C},S)$ 
is said to be \textbf{$\boldsymbol{\delta}$-2-transitive}  
if 
whenever 
$(c_1,c_2)$, and $(c_1',c_2')$ 
are
ordered pairs of chambers 
of $(\mathcal{C},S)$  
with
$\delta(c_1,c_2)=\delta(c_1',c_2')$, 
then the diagonal action of some element 
of $G$ 
maps  $(c_1,c_2)$ to $(c_1',c_2')$. 
\end{definition} 

If the action of $G$ on a building 
is $\delta$-2-transitive, 
then so is 
its restriction 
to the finite index subgroup of type-preserving automorphisms in $G$. 

\subsection{Consequences of $\delta$-2-transitivity}\label{subsec:delta-2-transitive=>}
In this subsection 
our main result is 
Theorem~\ref{thm:Main-delta-2-transitive}, 
which compares
the Davis-realization of a locally finite building 
with 
the space of compact open subgroups of 
a closed subgroup of its automorphism group 
acting $\delta$-2-transitively. 

The following proposition 
allows us 
to compute the distance between 
stabilizers of chambers 
for such groups. 

\begin{proposition}\label{prop:delta-2-transitive=>distance-formula(chamber_stab)}
Suppose that 
the action of a group $G$ 
on a locally finite building $(\mathcal{C},S)$  
with Weyl-distance $\delta$ 
is $\delta$-2-transitive and type-preserving. 
For each type $s$ 
let $q_s+1$ be the common cardinality of $s$-residues in $(\mathcal{C},S)$.  
Let $(c,c')$ be a pair of chambers of $(\mathcal{C},S)$ 
and let $s_1\cdots s_l$ be the type 
of some minimal gallery connecting $c$ to $c'$. 
Then $|G_c\colon G_c\cap G_{c'}|=\prod_{j=1}^l q_{s_j}$. 
\end{proposition}

\pf{}
Set $w=\delta(c,c')$. 
The index $|G_c\colon G_c\cap G_{c'}|$ 
is the length of the orbit of the chamber $c'$ 
under the action of $G_c$. 
Since the action of $G$ is type-preserving, 
this orbit is contained in the set 
$c_w=\{d\colon \delta(c,d)=\delta(c,c')\}$. 
Since the action of $G$ is $\delta$-2-transitive as well, 
the orbit is equal to this set. 
It remains to show that 
the cardinality of $c_w$ 
is equal to $\prod_{j=1}^l q_{s_j}$. 

Pick $d\in c_w$. 
By definition of buildings in terms of $W$-distance and by \cite[(3.1)v]{PiM7}, 
the chamber $c$ is connected to $d$ 
by a unique minimal gallery of type $s_1\cdots s_l$.  
On the other hand, any endpoint of a gallery of type $s_1\cdots s_l$ 
starting at $c$ 
will belong to $c_w$. 
Therefore 
the cardinality of $c_w$, 
hence the index $|G_c\colon G_c\cap G_{c'}|$, 
equals the number of %minimal 
galleries of type $s_1\cdots s_l$ 
starting at $c$. 
That latter number is equal to $\prod_{j=1}^l q_{s_j}$, 
establishing the claim. 
\qed
\vspace{3mm}

Before we continue  
with the preparation 
of Theorem~\ref{thm:Main-delta-2-transitive},   
we derive the following corollary. 

\begin{corollary}\label{rem:(ST)+type-preserving=>unimod}
A closed subgroup $G$ 
of the automorphism group of a locally finite building 
satisfying the conditions of  Proposition~\ref{prop:delta-2-transitive=>distance-formula(chamber_stab)}, 
is unimodular. 
In particular, the scale function of $G$ 
assumes the same value 
at a group element and its inverse.  
\end{corollary}
\pf{}
Let $\alpha$ be a bicontinuous automorphism 
of a totally disconnected, locally compact group $G$.  
Then the module of $\alpha$ 
equals $|\alpha(V)\colon \alpha(V)\cap V|\cdot |V\colon \alpha(V)\cap V|^{-1}$, 
where $V$ is an arbitrary compact, open subgroup of $G$. 

If we apply this observation 
in the situation of Proposition~\ref{prop:delta-2-transitive=>distance-formula(chamber_stab)} 
with $V=G_c$ 
and 
$\alpha$ being inner conjugation by $x\in G$, 
we conclude that 
the modular function of $G$ 
takes the value 
$ |G_{x.c}\colon G_{x.c}\cap G_{c}|\cdot |G_c\colon G_{x.c}\cap G_{c}|^{-1}$  
at $x$.  
Since traversing 
a minimal gallery joining $x.c$ to $c$ 
in the opposite order 
gives a minimal gallery joining $c$ to $x.c$, 
the formula for the index 
derived in  Proposition~\ref{prop:delta-2-transitive=>distance-formula(chamber_stab)}  
shows that 
$x$ has module $1$. 
Since $x$ was arbitrary, $G$ is unimodular. 

Our second claim 
follows from the first one,  
because 
the value of the modular function 
of $G$ 
at an automorphism $\alpha$ 
is $s_G(\alpha)/s_G(\alpha^{-1})$. 
We recall the proof. 
Choose $V$ to be tidy 
for the automorphism $\alpha$ 
of $G$. 
Then 
%the number 
$|\alpha(V)\colon \alpha(V)\cap V|$ equals the scale of $\alpha$ 
while 
%the number 
$|V\colon \alpha(V)\cap V|$ equals the scale of $\alpha^{-1}$. 
By the definition of the modular function $\Delta_G$ of $G$, 
we have $\Delta_G(\alpha)=s_G(\alpha)/s_G(\alpha^{-1})$ 
as claimed. 
\qed
\vspace{3mm} 

The following concept will be used to make the comparison in Theorem~\ref{thm:Main-delta-2-transitive}. 

\begin{definition}[adjacency graph of a chamber system]
Let $(\mathcal{C},S)$ be a chamber system 
and 
let $d : s \mapsto d_s$ be a map from $S$ to the positive real numbers. 
The \emph{adjacency graph of $(\mathcal{C},S)$ 
with respect to $(d_s)_{s \in S}$} 
is defined as follows. 
It is the labeled graph $\Gamma(\mathcal{C},S,(d_s)_{s \in S})$ 
having $\mathcal{C}$ as set of vertices; 
two vertices $c,c'\!\in\!\mathcal{C}$ are connected by an edge 
of label $s$ 
if and only if 
$c$ and $c'$ are $s$-adjacent for some $s\!\in\!S$; 
each edge  
of label $s$ 
is defined to have 
length $d_s$. 
\end{definition}

We are ready to state and prove our comparison theorem. 

\begin{theorem}[Comparison Theorem]
\label{thm:Main-delta-2-transitive}\label{prop:(ST)=> nice action}
Suppose that 
$(\mathcal{C},S)$ is a locally finite building 
with Weyl group distance function $\delta$ 
and let $\epsilon\in\{0,-1\}$. 
\begin{enumerate}
\item 
Suppose that $(d_s)_{s \in S}$ is a set of positive real numbers 
indexed by $S$. 
Then  
any map $\psi\colon |\mathcal{C}|_\epsilon \to \Gamma(\mathcal{C},S,(d_s)_{s \in S})$
which sends 
a point $x$ 
to some chamber $c\in \mathcal{C}$ 
such that 
$x\in|c|$ 
is a quasi-isometry. 
\item 
Let $G$ be a closed subgroup 
of the group of automorphisms of $(\mathcal{C},S)$ 
such that 
the action of $G$ 
is $\delta$-2-transitive. 
Then 
the map $\varphi\colon |\mathcal{C}|_\epsilon\to \bs G.$ 
mapping a point to its stabilizer is a quasi-isometric embedding. 
\end{enumerate}
\end{theorem}

\pf{}
We begin by proving the first claim. 
%
%The quantities 
Both 
$m:=\min\{d_s\colon s\in S\}$ 
and 
$M:=\max\{d_s\colon s\in S\}$ %, 
are finite 
and positive,  
because 
$S$ is a finite set 
and 
$\{d_s\colon s\in S\}$ 
is a set of positive numbers. 
The image of $\psi$ 
is obviously  
$M/2$-quasi-dense 
in $\Gamma(\mathcal{C},S,(d_s)_{s \in S})$ 
and 
we need to prove that 
$\psi$ is a quasi-isometric embedding 
as well. 

To see this, 
we compare distances 
between points 
in 
the space $|\mathcal{C}|_\epsilon$ 
and 
the graph $ \Gamma(\mathcal{C},S,(d_s)_{s \in S})$ 
to the gallery-distance between 
chambers corresponding to these points.  
To that end, 
denote by 
$d_\mathcal{C}$, $d_\Gamma$ and $d_\epsilon$ 
the gallery-distance on the set of chambers, 
the distance in the graph  $\Gamma(\mathcal{C},S,(d_s)_{s \in S})$  
and 
the distance on $|\mathcal{C}|_\epsilon$ 
respectively. 
Factor $\psi$ as 
the composite of a map 
$\psi'\colon (|\mathcal{C}|_\epsilon,d_\epsilon) \to (\mathcal{C},d_\mathcal{C})$,  
sending each point 
to some chamber containing it 
and 
the map $\iota\colon 
\bigl(\mathcal{C},d_\mathcal{C}\bigr) \to \bigl(\Gamma(\mathcal{C},S,(d_s)_{s \in S}),d_\Gamma\bigr)$ 
induced by 
the identity on $\mathcal{C}$. 
For any pair of chambers, 
$c$ and $c'$ 
we have 
\[
m\,d_\mathcal{C}(c,c')\le d_\Gamma(c,c')\le M\,d_\mathcal{C}(c,c')\,,
\]
which shows that $\iota$ is a quasi-isometric embedding. 

Furthermore, 
given two points,  
$x$ and $y$ say, 
in $|\mathcal{C}_\epsilon|$, 
there is a minimal gallery 
such that 
the geometric realizations 
of the chambers in that gallery 
do cover the geodesic 
joining $x$ to $y$. 
Denoting by $D$ 
the diameter of 
the geometric realization $|c|$ 
of a chamber 
in $|\mathcal{C}|_\epsilon$ 
and 
by $r$ 
the maximal gallery-distance 
between chambers 
in the same spherical residue 
we conclude that 
\[
d_\epsilon(x,y)-2rD\le D\,d_\mathcal{C}(\psi'(x),\psi'(y))\le d_\epsilon(x,y)+2rD\,. 
\]
It follows that 
$\psi'$ is a quasi-isometric embedding,  
and so 
$\psi$ is as well, 
finishing the proof 
of the first claim. 

%	The first claim follows from the following facts. 
%	First, there is a constant $m$ such that 
%	the geometric realization $|c|$ 
%	of any chamber has diameter $m$. 
%	Second, 
%	$S$ is finite so the set of numbers $\{ d_s \}_{s \in S}$ 
%	is bounded above and away from $0$. 

The second claim 
will be derived from the first 
%using 
and Proposition~\ref{prop:delta-2-transitive=>distance-formula(chamber_stab)}, 
once we have shown that 
we can 
reduce to the case where 
$G$ acts 
by type-preserving automorphisms. 
Let $G^\circ$ be the subgroup of 
type-preserving automorphisms in $G$. 
It is a closed subgroup of finite index, say $n$, in $G$. 
Therefore it is open in $G$, and $d(O,G^\circ\cap O)\leq \log(|G\colon G^\circ|)=\log n$ 
for each open subgroup $O$ of $G$. 
%Hence the map 
It follows that the map 
$\bs G^\circ.\to \bs G.$ 
induced by the inclusion $G^\circ\hookrightarrow G$ 
is a quasi-isometry. 
Since $G^\circ$ is $\delta$-2-transitive as well, 
we may --- and shall --- assume that $G=G^\circ$. 

Choose some map $\psi$ 
satisfying the conditions 
on the map with the same name 
in part~1. 
As before 
denote 
the flag 
consisting of 
the single vertex $\varnothing$ 
in $X$ 
by $(\varnothing)$. 
It defines a vertex of the simplicial complex underlying 
the model space $X$ of 
the Davis- and Moussong realizations. 
Denote 
the equivalence class  
with respect to the relation $\sim$ 
containing the pair $(c,x)\in \mathcal{C}\times X$ 
by $[c,x]$. 
Let 
$\mathcal{C}_{\varnothing}:=\{[c,|\varnothing|]\colon c\in\mathcal{C}\}$, 
a set of vertices 
of the simplicial complex underlying $|\mathcal{C}|_\epsilon$. 
%which was used 
%in the proof 
%of Proposition~\ref{prop:Aut(chamber-sys)=Isom(geom-realization)}. 
%
Let $q_s$ be the common cardinality of $s$-residues in $(\mathcal{C},S)$. 
Proposition~\ref{prop:delta-2-transitive=>distance-formula(chamber_stab)} implies that 
for the choice $d_s= 2\log q_s$ for $s\!\in\!S$ 
the map $\nu\colon \im(\varphi) \to \Gamma(\mathcal{C},S,(d_s)_{s \in S})$
defined by $\nu(G_{[c,|\varnothing|]}):= \psi([c,|\varnothing|])$ 
is an isometric embedding. 
Since the composite of 
the restriction of $\varphi$ to $\mathcal{C}_{\varnothing}$
with $\nu$ 
equals 
the restriction of $\psi$ to $\mathcal{C}_{\varnothing}$
it follows that 
the restriction of $\varphi$ to $\mathcal{C}_{\varnothing}$ 
is a quasi-isometric embedding, 
because $\psi$ is 
by part~1, 
which we already proved. 

It follows that 
$\varphi$ is a quasi-isometric embedding as well 
because 
$\mathcal{C}_{\varnothing}$ is quasi-dense in $|\mathcal{C}|_\epsilon$ 
and 
the distance between 
the stabilizer of a point %$x$ 
in the geometric realization $|c|$ 
of a chamber $c$ 
and the stabilizer of the point $[c,|\varnothing|]\in \mathcal{C}_{\varnothing}\cap |c|$ 
is bounded above by a constant $M$ 
independent of $c$ 
since  
$|c|$ is a finite complex 
and that $G$ acts transitively on chambers. 
\qed 

\subsection{Consequences of the Comparison Theorem}\label{subsec:ComparisonThm=>}
The most important consequence of the Comparison Theorem, 
Theorem~\ref{thm:Main-delta-2-transitive}, 
is the inequality between the flat rank  of a $\delta$-2-transitive automorphism group 
and the rank of the building the group acts on. 
This result follows from the following

\begin{proposition}\label{prop:nice-action=>flat-subgr->quasi-flat}
Let $G$ be a totally disconnected, locally compact group. 
Suppose that $G$ acts 
on a metric space $X$ 
in such a way that 
$G$-stabilizers 
of points are compact, open subgroups 
of $G$.  
Assume that 
the map $X\to \bs G.$, 
which assigns a point its stabilizer, 
is a quasi-isometric embedding. 
Let $H$ be a flat subgroup 
of $G$ 
of finite flat rank  $n$. 
Then, 
for any point $x$ in $X$ 
the inclusion of 
the $H$-orbit of $x$ defines an $n$-quasi-flat in $X$. 
\end{proposition}

\pf{}
Let $x$ be any point in $X$. 
Since the map $x\mapsto G_x$ is an equivariant quasi-isometric embedding, 
the orbit of $x$ under $H$ is quasi-isometric to 
the orbit of its stabilizer $G_x$ under $H$ acting by conjugation. 
The latter orbit is quasi-isometric to the $H$-orbit 
of a tidy subgroup, say $O$, for $H$. 
But $H.O$ is isometric to the subset $\mathbb{Z}^n$ 
of  $\mathbb{R}^n$ with the norm $\|\cdot\|_H$ 
introduced in Subsection~\ref{subsec:fundamentals-tidy}. 
The subset  $\mathbb{Z}^n$ is quasi-dense in  
$\mathbb{R}^n$ equipped with that norm 
and therefore, 
we obtain a quasi-isometric embedding of 
$\mathbb{R}^n$ equipped with $\|\cdot\|_H$ into $X$. 
But the identity map between 
$\mathbb{R}^n$ equipped with $\|\cdot\|_H$
and 
$\mathbb{R}^n$ with the Euclidean norm 
is bi-Lipschitz. 
Composing with this map, we obtain 
an $n$-quasi-flat in $X$ 
as claimed. 
\qed
\vspace{3mm} 

Before deriving the rank inequality, we prove the following characterization for existence of fixed points. 

\begin{corollary}\label{cor:nice-action=> elliptic are in G(1)}
Let $G$ be a totally disconnected, locally compact group. 
Suppose that 
$G$ acts isometrically 
on a complete $\CAT0$--space $X$ 
with compact, open point stabilizers. 
Assume that 
the map $X\to \bs G.$, 
which assigns a point its stabilizer, 
is a quasi-isometric embedding. 
Then 
an element $g$ of $G$ 
has a fixed point in $X$ 
if, and only if, $s_G(g)=1=s_G(g^{-1})$. 
\end{corollary}

\pf{}
Let $g$ be an element of $G$
and 
let $x$ be a point of $X$. 
The subgroup $\langle g\rangle$ 
is a flat subgroup of $G$ 
of flat rank  $0$ or~$1$ 
and we have 
$\frk(\langle g\rangle)=0$ 
if and only if 
$s_G(g)=1=s_G(g^{-1})$. 
Hence 
to prove our claim, 
we need to show that 
$g$ has a fixed point in $X$ 
if and only if 
$\frk(\langle g\rangle)=0$. 

Proposition~\ref{prop:nice-action=>flat-subgr->quasi-flat} 
is applicable 
with $H$ equal to $\langle g\rangle$. 
Therefore 
the set $\langle g\rangle.x$ 
is quasi-isometric to a point or the real line
in the cases $\frk(\langle g\rangle)=0$ and $\frk(\langle g\rangle)=1$ 
respectively. 
We conclude that 
the set $\langle g\rangle.x$ 
is bounded 
if and only if 
 $\frk(\langle g\rangle)=0$. 
 
Since $\langle g\rangle$ acts by isometries 
on the complete $\CAT0$-space $X$, 
if it has a bounded orbit, 
it admits a fixed point by \cite[II.2, Corollary~2.8(1)]{GMW319}. 
The converse 
of the latter statement 
is trivial. 
We conclude that 
$g$ has a fixed point in $X$ 
if and only if 
$\frk(\langle g\rangle)=0$ 
as had to be shown.  
\qed
\vspace{3mm} 

We adopt the following definition for the rank of a complete $\CAT0$-space. 
For alternative definitions see \cite[pp.~127-133]{asymp-inv(infG)}. 

\begin{definition}
The \emph{rank} of a complete $\CAT0$-space $X$, denoted $\rk(X)$, is the maximal dimension of an isometrically embedded Euclidean space in $X$. 
\end{definition}

Recall that 
a metric space $X$ 
is called \emph{cocompact} 
if and only if 
the isometry group of $X$ 
acts cocompactly on $X$ \cite[p. 202]{GMW319}. 

\begin{theorem}
\label{thm:frk(G) at most rk(X)}
Let $G$ be a totally disconnected, locally compact group. 
Suppose that $G$ acts isometrically  
on a complete, locally compact, cocompact $\CAT0$-space $X$ 
with compact, open point stabilizers. 
Assume that the map $X\to \bs G.$, 
which assigns a point its stabilizer, 
is a quasi-isometric embedding. 
Then 
$\frk(G)\leq \rk(X)$; 
in particular 
the flat rank  of $G$ is finite. 
We have 
$\frk(G)=0$ 
if, and only if, 
every element of $G$ fixes a point in $X$. 
\end{theorem} 

\pf{}
The hypotheses on $X$ guarantee that 
the rank of $X$ is finite 
and equals the maximal dimension of quasi-flats in $X$ 
by \cite[Theorem~C]{local-struc(lengthS_curv<)}.  
The hypotheses on the action of $G$ on $X$ 
enable us to apply Proposition~\ref{prop:nice-action=>flat-subgr->quasi-flat},  
which, 
together with the first observation of this proof, 
implies that 
$\frk(G)\leq \rk(X)<\infty$. 
The last statement follows from Corollary~\ref{cor:nice-action=> elliptic are in G(1)}. 
\qed
\vspace{3mm} 

The rank of a Gromov-hyperbolic $\CAT0$-space is $1$. 
This leads to the following special case of Theorem \ref{thm:frk(G) at most rk(X)}. 

\begin{corollary}\label{cor:frk(G) for hyp-CAT0}
Let $G$, $X$ and $X\to \bs G.$ be as above. 
Assume further that $X$ is Gromov-hyperbolic. 
Then $\frk(G)=1$, unless every element of $G$ has a fixed point in $X$, in which case $\frk(G)=0$. 
\qed
\end{corollary} 
%	
%	This corollary covers most examples treated in \cite{G-flatrk_le1}. 

\subsection{Equality of the rank of a building and the rank of an apartment}
\label{subsec:rk(building)=rk(WeylG)}
The purpose of this subsection is to prove that 
the Davis-realization of a building 
has the same rank as any of its apartments. 
We do not claim that 
any flat of the building 
is contained in an apartment, 
though this is probably true as well. 
This stronger statement 
is known to be true 
for Euclidean buildings 
by Theorem~1 in \cite[Chapter~VI, Section~7]{buildings}. 
We believe,  
we can show it to be true also 
if the building is 
Moufang 
and locally finite.  

% I let you choose the final version of this subsection 
% if you change something, please change the introduction accordingly

\begin{proposition}\label{prop:rk(building)=rk(apartment)}
Let $(\mathcal{C},S)$ be a 
building 
with $S$ finite 
and 
Weyl group $W$. 
If $|\mathcal{C}|_0$ contains a $d$-flat, 
then so does $|W|_0$. 
Hence $\rk(|\mathcal{C}|_0)=\rk(|W|_0)$. 
\end{proposition}
\pf{}
Let $F$ be a $d$-flat in $|\mathcal{C}|_0$. 
Applying Lemma~9.34 in chapter~II of \cite{GMW319} 
with $Y$ equal to $\mathbb{R}^d$ 
and $X$ equal to $|W|_0$ 
and 
using the isomorphism of 
the geometric realization of 
any apartment with $|W|_0$, 
we see that 
it suffices to show that 
for each $n$ in $\mathbb{N}$ 
there is an apartment $A_n$ of $(\mathcal{C},S)$ 
which contains 
an isometric copy of 
the ball of radius $n$ around $\mathbf{0}$ in $\mathbb{R}^d$. 

This isometric copy, 
$B_n$, 
of the ball of radius $n$ around $\mathbf{0}$ in $\mathbb{R}^d$ 
will be taken to lie 
inside $F$. 
Let $n$ be a natural number, 
$o$ some point in $F$ 
and 
$B_n$ the ball of radius $n$ around $o$ in $F$.
To show that 
there is an apartment $A_n$ 
containing $B_n$, 
we will prove that 
there are two chambers $c_n$ and $c'_n$ 
such that 
the minimal galleries 
connecting $c_n$ and $c'_n$ 
cover $B_n$. 
Then by combinatorial convexity any apartment $A_n$ 
containing $c_n$ and $c'_n$ 
contains $B_n$ 
and 
by our introductory remark 
the proposition 
follows, 
because $n$ was arbitrary. 

To determine 
how we should choose 
the chambers $c_n$ and $c'_n$, 
we first 
take a look at 
the way walls in $|\mathcal{C}|_0$ 
intersect the flat $F$. 
Since 
any geodesic 
joining two points 
of a wall 
lies entirely inside 
that wall, 
% walls are defined as fixed point sets of reflections in apartments; 
% these reflections act by isometries and geodesics are unique 
the intersection of a wall with $F$ 
is an affine subspace of $F$. 
(It can be shown that 
the affine subspaces 
of $F$ 
arising in this way 
are either empty 
or of codimension at most $1$ in $F$, 
but we will not make use 
of this additional information.) 
Note further that 
the family 
of affine subspaces 
arising as 
intersections of 
walls in $|\mathcal{C}|_0$ 
with 
$F$ 
is locally finite. 

If 
$M$ is a wall in $|\mathcal{C}|_0$ 
and 
two points $p$ and $p'$ in $F$ 
are separated by $M\cap F$ 
in $F$,  
then 
$p$ and $p'$ 
are separated by $M$ in $|\mathcal{C}|_0$. 
Therefore, 
if we choose 
chambers $c_n$ and $c'_n$ 
to contain 
points $p_n$ and $p'_n$ 
in $F$ 
such that 
no intersection of a wall with $F$ 
separates any point in $B_n$ 
from both $p_n$ and $p'_n$, 
then 
the minimal galleries 
connecting $c_n$ and $c'_n$ 
cover $B_n$. 
The following lemma 
demonstrates that 
such a choice 
of $p_n$ and $p'_n$ 
is always possible. 
This concludes 
the proof 
modulo Lemma~\ref{lem:pos-observation_rel-hyperplanes}. 
\qed 

The following Lemma completes the proof of Proposition~\ref{prop:rk(building)=rk(apartment)}. 

\begin{lemma}\label{lem:pos-observation_rel-hyperplanes}
Let $B$ be a bounded, convex 
subset of  $\mathbb{R}^d$ 
with non-empty interior 
and 
let $\mathcal{M}$ be a locally finite 
collection of hyperplanes in $\mathbb{R}^d$. 
Then there exist two points $p$ and $p'$ 
in the complement of $\bigcup \mathcal{M}$ 
such that 
no element of $\mathcal{M}$ 
separates 
any point in $B$  
from both $p$ and $p'$. 
\end{lemma}
\pf{}
Let $\mathcal{M}_\cap$ be 
the subfamily of $\mathcal{M}$ 
consisting of those members of $\mathcal{M}$ 
which intersect the interior of $B$. 
Since $\mathcal{M}$ is locally finite 
and $B$ is bounded, 
$\mathcal{M}_\cap$ is finite. 

We will choose 
the points $p$ and $p'$ 
from the complement of $\bigcup \mathcal{M}$ 
to have 
the following two properties. 
\begin{enumerate}
\item 
The line segment $[p,p']$ intersects the interior of $B$. 
\item 
Each element of $\mathcal{M}_\cap$ 
separates $p$ from $p'$. 
\end{enumerate}
The first property 
ensures that 
for every element $M$ of $\mathcal{M}\setminus \mathcal{M}_\cap$ 
at least one of the points $p$ and $p'$ 
lies on the same side of $M$ as $B$,  
while the second property 
ensures that 
for every element $M$ of $\mathcal{M}_\cap$ 
any point of $B$ 
lies 
in the same closed halfspace 
with respect to $M$ 
as either $p$ or $p'$. 
Therefore 
the two conditions ensure 
that 
no element of $\mathcal{M}$ 
separates 
any point in $B$  
from both $p$ and $p'$, 
and it will suffice to conform with 
the conditions 1 and~2 
to conclude the proof. 

Since $\mathcal{M}_\cap$ is finite, 
there is at least one unbounded component, 
say $C$, 
of the complement of $\bigcup\mathcal{M}_\cap$ 
which contains interior points of $B$. 
This can be shown by induction on the number of elements in $\mathcal{M}_\cap$. 
Let $p$ be 
any interior point of $B$  
inside $C$. 
By the definition of $\mathcal{M}_\cap$, 
the point $p$ 
does not lie in any element of $\mathcal{M}$. 
Note that 
% whenever instead of however, I guess
whenever the point $p'$ 
is subsequently chosen, 
the pair $(p,p')$ 
will satisfy property~1 
by our choice of $p$. 

Let $H$ be the collection of all open halfspaces 
with respect to $\mathcal{M}_\cap$, 
which do not contain $p$. 
Since the component $C$ 
of the complement of $\bigcup\mathcal{M}_\cap$ 
which contains $p$ 
is unbounded, 
$\bigcap H$ is not empty. 
The set  $\bigcap H$ 
is also open in $\mathbb{R}^d$. 
Since $\mathcal{M}$ is locally finite, 
we can therefore choose a point $p'$ 
in $\bigcap H$ not contained in $\bigcup \mathcal{M}$. 
By construction, 
the pair $(p,p')$ 
satisfies property~2 
and 
we have already noted that 
it satisfies property~1 
as well.  
As seen above, 
this proves that $p$ and $p'$ 
have the property sought 
in the statement of the Lemma. 
\qed

% SECTION 3

\section{Algebraic rank as a lower bound}\label{sec:lower_bound}
In this section, we show that when $G$ is either 
a semi-simple algebraic group over a local field  
or 
a topological Kac-Moody group over a finite field, 
the algebraic rank of the Weyl group $W$ (Definition~\ref{def:algrk(G)}) 
is a lower bound for the flat rank  of $G$. 

One strategy to get a lower bound of $\frk(G)$ of geometric nature, i.e. coming from the building $X$, is to use the stabilizer map 
%$\varphi:x\mapsto {\rm Stab}_G(x)$. 
$\varphi:x\mapsto G_x$. 
Indeed, according to Theorem~\ref{thm:Main-delta-2-transitive}, $\varphi$ maps an $n$-quasi-flat in $|\mathcal{C}|_\epsilon$ to an $n$-quasi-flat in $\bs G.$, which can even be assumed to consist of vertex stabilizers. 
Still, in order to make this strategy work, one needs a connection between flats in the space $\bs G.$ and flat subgroups in $G$ itself. 
One connection is stated by the following 

\begin{conjecture}
Let $\mathcal{H}$ be a group of automorphisms 
of a totally disconnected, locally compact group $G$.  
Suppose that some (equivalently every) orbit of $\mathcal{H}$ 
in $\bs G.$ is quasi-isometric to $\mathbb{R}^n$. 
Then $\mathcal{H}$ is flat, of flat rank  $n$. 
\end{conjecture}

This conjecture has been verified for $n=0$ \cite[Proposition~5]{direction(aut(tdlcG))}. 
Going back to the announced algebraic lower bound, we now introduce the algebraic rank of a group. 

\begin{definition}\label{def:algrk(G)}
Let $H$ be a group. 
The algebraic rank of $H$, denoted $\ark(H)$, is the supremum of the ranks of the free abelian subgroups of $H$. 
\end{definition}

In order to achieve the aim of this section, we will show that given a free abelian subgroup $A$ of the Weyl group $W$ of the building $X$, we can lift $A$ to a flat subgroup of the isometry group $G$. 
This holds in both classes of groups considered. 

\begin{theorem}\label{thm:lifting-free_ab<(W)}
\rule{0pt}{0pt}
\begin{enumerate}
\item 
Let $\mathbf{G}$ be an algebraic semi-simple group $\mathbf{G}$ over a local field $k$, 
and let $G=\mathbf{G}(k)$ be its rational points. 
Let $\mathbf{S}$ be a maximal $k$-split torus in $\mathbf{G}$, 
and let $W=W_{{\rm aff}}(\mathbf{S},\mathbf{G})$ be the affine Weyl group of $\mathbf{G}$ 
with respect to $\mathbf{S}$.  
Then $\mathbf{S}(k)$ is a flat subgroup of $G$, of flat rank $\ark(W)= k$-$\rk(\mathbf{G})$. \\
In particular, we have: $\ark(W)\le \frk(G)$. 
\item 
Let $G$ be a group with a locally finite twin root datum 
of associated Weyl group $W$. 
Let $A$ be an abelian subgroup of $W$ 
and let $\widetilde{A}$ be the inverse image of $A$ 
under the natural map $N\to W$. 
Then 
$\widetilde{A}$ is a flat subgroup of $\overline{G}$ 
and 
$\frk(\widetilde{A})=\rank_\mathbb{Z}(A)$. \\
In particular, we have: $\ark(W)\le \frk(G)$. \end{enumerate}
\end{theorem}

Before beginning the proof of the above theorem, we note the following corollary, which was obtained by different means in unpublished work of the first and third authors. 

\begin{corollary}
\label{cor:flatrk(semi-simpleG)}
Let $\mathbf{G}$ be an algebraic semi-simple group over a local field $k$, with affine Weyl group $W$. 
Then $\frk(\mathbf{G}(k))=k$-$\rk(\mathbf{G})=\ark(W)$. 
\end{corollary}

This result and its proof are valid for any closed subgroup lying between $\mathbf{G}(k)$ and its closed subgroup $\mathbf{G}(k)^+$ generated by unipotent radicals of parabolic $k$-subgroups 
--- see e.g. \cite[I.2.3]{EdM3.17} for a summary about $\mathbf{G}(k)^+$. 
It suffices to replace $\mathbf{S}(k)$ by $\mathbf{S} \cap \mathbf{G}(k)^+$. 

\vspace{3mm} 

\pf{}
Let 
$\mathbf{S}$ be a maximal $k$-split torus of $\mathbf{G}$ 
and 
let $W$ be the affine Weyl group of $\mathbf{G}$ with respect to $\mathbf{S}$. 
Part~1 of Theorem~\ref{thm:lifting-free_ab<(W)} 
shows that \linebreak 
$k$-$\rk(\mathbf{G})=\ark(W)\le \frk(\mathbf{G}(k))$. 

To show the inequality $\ark(W)\ge \frk(\mathbf{G}(k))$, 
note that 
the action of $\mathbf{G}(k)$ 
on its affine building $X$ 
is given by a BN-pair in $\mathbf{G}(k)$, 
which implies that 
this action is strongly transitive.  
In particular,  
the action of $\mathbf{G}(k)$ on $X$ 
is $\delta$-2-transitive 
for the canonical $W$-metric $\delta$. 
Part~2 of Theorem~\ref{thm:Main-delta-2-transitive} shows that 
the map 
assigning a point of $X$ its stabilizer in $\mathbf{G}(k)$ 
is a quasi-isometric embedding. 
Hence 
Theorem~\ref{thm:frk(G) at most rk(X)} is applicable 
and yields 
$\frk(\mathbf{G}(k)) \le \rk(X)$. 
In the case at hand, we have $\rk(X)=\ark(W)$ 
because maximal flats in $X$ are apartments in $X$ \cite[VI.7]{buildings}, 
and $W$ is virtually free abelian. 

We finally conclude that 
$k$-$\rk(\mathbf{G})=\ark(W)= \frk(\mathbf{G}(k))$. 
\qed
\vspace{3mm}

\pf{Theorem~\ref{thm:lifting-free_ab<(W)}}
We treat separately cases~1 and~2. 

\vspace{3mm}

In case~1, the group $\mathbf{S}(k)$ is flat by Theorem~5.9 of \cite{tidy<:commAut(tdlcG)}, being topologically isomorphic to a power of the multiplicative group of $k$, which itself is generated by the group of units and a uniformizer of $k$. 

It remains to show that $\frk(\mathbf{S}(k))=\rank_\mathbb{Z}\left(\mathbf{S}(k)/(\mathbf{S}(k))(1)\right)$ 
equals $\ark(W)=k$-$\rk(\mathbf{G})$. 
As seen in the proof of Corollary~\ref{cor:flatrk(semi-simpleG)}, 
the map 
assigning 
to a point of the affine building $X$ 
its stabilizer in $\mathbf{G}(k)$ 
is a quasi-isometric embedding.
Hence 
Corollary~\ref{cor:nice-action=> elliptic are in G(1)} 
can be applied to 
the action of $G=\mathbf{G}(k)$ on $X$. 
Therefore 
$(\mathbf{S}(k))(1)$ is the subgroup of $\mathbf{S}(k)$ 
of elements admitting fixed points in $X$. 
Denote by $A_\mathbf{S}$ 
the affine apartment of $X$ 
associated to $\mathbf{S}$/ 
The group $\mathbf{S}(k)$ 
leaves $A_\mathbf{S}$ invariant. 
Since $A_\mathbf{S}$ is a convex subspace of the $\CAT0$-space $X$, every element of $\mathbf{S}(k)$ 
admitting a fixed point in $X$ fixes a point of $A_\mathbf{S}$ 
(the projection of the fixed point onto $A_\mathbf{S}$). 
Therefore 
$\mathbf{S}(k)/(\mathbf{S}(k))(1)$ identifies with 
the group of translations of $A_\mathbf{S}$ induced by $\mathbf{S}(k)$, 
which is a subgroup of finite index in the translation lattice of $W$, 
which is itself an abelian subgroup of $W$ of maximal rank. 
Therefore 
$\frk(\mathbf{S}(k))=\ark(W)$, 
which coincides with $k$-$\rk(\mathbf{G})$.  
This settles case~1. 

\vspace{3mm}

We now reduce case~2 
to Proposition~\ref{prop:almost abelian=>flat}, 
which we will prove later. 
Let $(H,(U_a)_{a\in\Phi})$ be the locally finite root datum in $G$,  
$(\mathcal{C},S)$ the positive building of $(H,(U_a)_{a\in\Phi})$ 
and 
$\mathcal{A}$ the standard apartment in $(\mathcal{C},S)$.
Let $\delta$ be the $W$-distance on $(\mathcal{C},S)$ 
and 
let $\epsilon\in\{0,-1\}$. 

Since the group $H$ is finite,  
$\widetilde{A}$ is finitely generated. 
Since $A$ is abelian,   
the commutator of two elements of $\widetilde{A}$ 
is contained in $H$, 
hence is of finite order. 
It follows that 
$\widetilde{A}$ is a flat subgroup of $\overline{G}$ 
by Proposition~\ref{prop:almost abelian=>flat}. 

The $G$-action on $(\mathcal{C},S)$ is $\delta$-2-transitive, 
hence so is the $\overline{G}$-action.  
Hence the map 
sending a point in $|\mathcal{C}|_\epsilon$ to its stabilizer in $\overline{G}$ 
is a quasi-isometric embedding 
by Theorem~\ref{prop:(ST)=> nice action}.
Corollary~\ref{cor:nice-action=> elliptic are in G(1)} 
then implies that 
the subgroup $\widetilde{A}(1)$ 
is the set of elements in $\widetilde{A}$ which fix some point in $|\mathcal{C}|_\epsilon$. 

The group $N$ leaves $\mathcal{A}$ invariant and the induced action of $N$ on $|\mathcal{A}|_\epsilon$ is equivariant to the action of $W$ when $\mathcal{A}$ is identified with the Coxeter complex of $W$. 
Since the kernel of the action of $N$ on $|\mathcal{A}|_\epsilon$ is the finite group $H$, the action of $N$ on $|\mathcal{A}|_\epsilon$ is proper. 
Hence an element of infinite order in $N$ does not have a fixed point in $|\mathcal{A}|_\epsilon$ and therefore has no fixed point in $|\mathcal{C}|_\epsilon$ either, because $|\mathcal{A}|_\epsilon$ is an $N$-stable convex subspace of the $\CAT0$-space $|\mathcal{C}|_\epsilon$. 
Therefore $\widetilde{A}(1)$ is the set of elements of $\widetilde{A}$ of finite order. 
This implies that 
$\rank_\mathbb{Z}(\widetilde{A}/\widetilde{A}(1))=\rank_\mathbb{Z}(A)$, as claimed. 
\qed
\vspace{3mm}

The remainder of this section is devoted  to the proof of the following

\begin{proposition}\label{prop:almost abelian=>flat}
Let $G$ be a totally disconnected locally compact group 
and 
let $A$ be a subgroup of $G$ 
which is a finite extension of a finitely generated abelian group. 
Then $A$ is a flat subgroup of $G$. 
\end{proposition}

We split the proof into several subclaims. 

\begin{lemma}\label{lem:cp=>flat}
If $C$ is a compact subgroup 
of $G$, 
then 
there is a base of neighborhoods of $e$ 
consisting of $C$-stable, compact, open subgroups. 
In particular, any compact subgroup is flat.
\end{lemma}

\pf{}
If $V$ is a compact, open subgroup of $G$, then $O=\bigcap_{c\in C} cVc^{-1}$ is a compact, open subgroup 
of $G$ which is contained in $V$ and is normalized by $C$. 
Since $G$ has a base of neighborhoods consisting of compact, open subgroups, this proves the first claim. 
The second claim follows from it.  
\qed 
\vspace{3mm} 

We remind the reader on the tidying procedure 
defined in~\cite{tidy<:commAut(tdlcG)}
which 
for any automorphism $\alpha$ of $G$ 
outputs a subgroup tidy for $\alpha$. 
It will be used in the proofs of 
Lemmas~\ref{lem;tidy_cp-is-aut-stable->cp-stable} 
to~\ref{lem:almost-abelianG has local tidy subgroups}.

\begin{algorithm}[$\alpha$-tidying procedure]\label{algo:tidy3}
~\\
{[0]} 
Choose a compact open subgroup $O\leqslant G$.\hfill\\
{[1]}
Let $^kO:=\bigcap_{i=0}^k \alpha^i(O)$.
For some $n\in \mathbb{N}$ (hence for all $n'\ge n$) we have \\
\hphantom{{[1]}}
$^nO=\bigl(\bigcap_{i\ge 0} \alpha^i({^nO})\bigr)\cdot \bigl(\bigcap_{i\ge 0} \alpha^{-i}({^nO})\bigr)$. 
Set $O'={^nO}$.\\
{[2]} 
For each compact, open subgroup $V$ let\\   
\hphantom{{[2]}}
$
\mathcal{K}_{\alpha,V}:= \bigl\{k\in G \colon \{\alpha^n(k)\}_{n\in\mathbb{Z}}\ 
\text{is bounded and $\alpha^n(k)\in V$ for all large $n$}\bigr\}
$. 
\hphantom{{[2]}}
Put $K_{\alpha,V} = \overline{\mathcal{K}}_{\alpha,V}$
and 
$K_\alpha =  \bigcap\{K_{\alpha,V} \colon V\  \text{is a compact, open subgroup}\}.$\\ 
{[3]} 
Form
$O^*:=\{x\in O'\colon kxk^{-1}\in O'K_\alpha\,\,\forall k\in K_\alpha\}$ 
and define $O'':=O^* K_\alpha$. \\
\hphantom{{[2]}}
The group $O''$ is tidy for $\alpha$ and we output $O''$.
\end{algorithm}

\begin{lemma}\label{lem;tidy_cp-is-aut-stable->cp-stable}
Let $\alpha$ be an automorphism of $G$ 
and 
$C$ an $\alpha$-stable compact subgroup of $G$.  
Choose 
a $C$-stable compact, open subgroup $O$
of $G$ as in Lemma~\ref{lem:cp=>flat}.  
Then 
$\alpha(O)$ 
and $K_\alpha$ 
are $C$-stable. 
Hence 
the output 
derived from 
applying the $\alpha$-tidying procedure 
in Algorithm~\ref{algo:tidy3} 
to $O$  
is $C$-stable. 
\end{lemma} 
\pf{}
To show that 
$\alpha(O)$ is $C$-stable, 
let $c$ be an element of $C$. 
Then 
\[
c\alpha(O)c^{-1}=\alpha\bigl(\alpha^{-1}(c)O\alpha^{-1}(c)^{-1}\bigr)=\alpha(O)\,,  
\]
which shows that 
$\alpha(O)$ is $C$-stable. 

Towards proving that 
$K_\alpha$ is $C$-stable, 
we first prove 
that 
if $C$ is $\alpha$-stable 
and 
$V$ is a $C$-stable compact, open subgroup 
of $G$ 
then 
$\overline{\mathcal{K}}_{\alpha,V}$ is $C$-stable. 
Let $c$ be an element of $C$ 
and 
$k$ an element of $\mathcal{K}_{\alpha,V}$. 
Then, 
for all $n\in\mathbb{Z}$, 
\[
\alpha^n(ckc^{-1})=
\alpha^n(c)\,\alpha^n(k)\,\alpha^n(c^{-1})\subseteq 
C\, \{\alpha^n(k)\colon n\in\mathbb{Z}\}\,C \,,
\] 
and hence 
the set $\{\alpha^n(ckc^{-1})\colon n\in\mathbb{Z}\}$ 
is bounded. 
Further, 
for all $n$ such that 
$\alpha^n(k)\in V$ 
we have 
\[
\alpha^n(ckc^{-1})\in \alpha^n(c)V\alpha^n(c)^{-1}=V
\]
proving that 
indeed 
$\overline{\mathcal{K}}_{\alpha,V}$ is $C$-stable 
if $C$ is $\alpha$-stable 
and 
$V$ is a $C$-stable. 

To derive that 
$K_\alpha$ is $C$ stable,  
note that 
whenever $V'\leqslant V$ 
are compact, open subgroups  
of $G$ 
then $\overline{\mathcal{K}}_{\alpha,V'}\leqslant \overline{\mathcal{K}}_{\alpha,V}$. 
Hence if $\mathcal{O}$ is a base of neighborhoods of $e$, 
consisting of compact, open subgroups 
of $G$ 
then 
$K_\alpha=\bigcap\{\overline{\mathcal{K}}_{\alpha,V}\colon V\in\mathcal{O} \}$. 
Since $C$ is compact, 
there is a base of neighborhoods $\mathcal{O}$ of $e$ 
consisting of $C$-stable, compact, open subgroups 
by Lemma~\ref{lem:cp=>flat}. 
We conclude that 
$K_\alpha$ 
is the intersection of 
a family of $C$-stable sets, 
hence is itself $C$-stable. 

Finally, 
we establish 
that the output of Algorithm~\ref{algo:tidy3} 
is $C$-stable. 
The group $O'$ 
is $C$-stable 
since 
it is the intersection of $C$-stable subgroups. 
Therefore 
$O^*$ is $C$-stable as well. 
Since $K_\alpha$ 
is $C$-stable, 
it follows that 
the output $O''$ 
is also $C$-stable 
as claimed. 
\qed 

\begin{lemma}\label{lem:ind-step aut-com-mod-cp}
Let $\mathcal{A}$ be a set of automorphisms of $G$, 
$C$ a compact subgroup of $G$ 
stable under each element of  $\mathcal{A}$ 
and let $\gamma$ be an automorphism of $G$ 
stabilizing $C$ 
such that 
$[\gamma,\alpha]$ is an inner automorphism in $C$ 
for each $\alpha$ in $\mathcal{A}$.  
Suppose there is a common tidy subgroup for $\mathcal{A}$ 
which is $C$-stable. 
Then there is a common tidy subgroup for $\mathcal{A}\cup\{\gamma\}$ 
which is $C$-stable. 
\end{lemma}
\pf{}
Let $O$ be a common tidy subgroup for $\mathcal{A}$ 
which is $C$-stable. 
We will show that 
the output $O''$ of 
the $\gamma$-tidying procedure, Algorithm~\ref{algo:tidy3},  
on the input $O$ 
produces a common tidy subgroup for $\mathcal{A}\cup\{\gamma\}$ 
which is $C$-stable. 
Lemma~\ref{lem;tidy_cp-is-aut-stable->cp-stable} 
shows that 
$O''$ is $C$-stable 
and we have to prove 
it is tidy for each element $\alpha$ in $\mathcal{A}$. 

Since $[\gamma,\alpha]\in C$, 
for all $C$-stable, 
compact, open subgroups $V$, 
we have  
\[
|\alpha\gamma(V)\colon \alpha\gamma(V)\cap \gamma(V)|=
|\gamma(\alpha(V))\colon \gamma(\alpha(V)\cap V)|=
|\alpha(V)\colon \alpha(V)\cap V|\,.
\]
Hence, 
if $V$ is $\alpha$-tidy and $C$-stable, 
then  
$\gamma(V)$ is $\alpha$-tidy 
and 
it is $C$-stable 
by Lemma~\ref{lem;tidy_cp-is-aut-stable->cp-stable} 
as well. 
Using this observation, 
induction on $i$ shows that 
$\gamma^i(O)$ 
is $\alpha$-tidy 
for each $i\in\mathbb{N}$. 
Since any finite intersection of 
$\alpha$-tidy subgroups 
is $\alpha$-tidy 
by Lemma~10 in \cite{tdlcG.structure}, 
the output $O'$ of step~1  
of Algorithm~\ref{algo:tidy3} 
will be $\alpha$-tidy. 

We will show next that 
$\alpha(K_\gamma)=K_\gamma$. 
Theorem~3.3 in~\cite{tidy<:commAut(tdlcG)} 
then implies that 
$O''$ is tidy for $\alpha$, 
finishing the proof. 

Towards proving 
our remaining claim $\alpha(K_\gamma)=K_\gamma$,  
we now show that 
if %whenever 
$V$ is a $C$-stable, compact, open subgroup of $G$ 
then  
$\alpha(\mathcal{K}_{\gamma,V})=\mathcal{K}_{\gamma,\alpha(V)}$. 
Using 
our assumptions 
$[\gamma,\alpha]\subseteq C$ 
and $\gamma(C)=C$, 
the equation 
\[
[\gamma^n,\alpha]=
[\gamma,\alpha]^{\gamma^{n-1}}\cdot [\gamma,\alpha]^{\gamma^{n-2}}\cdots 
[\gamma,\alpha]^{\gamma^{n-1}}\cdot [\gamma,\alpha]
\]
shows that 
for all 
%integers 
$n$ in $\mathbb{Z}$ 
there are 
%elements 
$c_n$ in $C$ 
such that 
$\gamma^n\alpha=\kappa(c_n)\,\alpha\gamma^n$,  
where $\kappa(g)$ denotes 
conjugation by 
%the element 
$g$. 
Therefore, 
if $k\in \mathcal{K}_{\gamma,V}$ 
%for some compact, open subgroup $V$ 
then 
\[
\gamma^n(\alpha(k))= 
c_{n}(\alpha\gamma^n(k))c_{n}^{-1}
\subseteq C \alpha(\gamma^n(k)) C
\] 
is bounded. 
Since $\alpha(C)=C$ 
we have %conclude that 
$\gamma^n\alpha=\kappa(c_n)\,\alpha\gamma^n=\alpha\kappa(\alpha^{-1}(c_n))\,\gamma^n$ 
for all $n$. 
Put $c_n'=\alpha^{-1}(c_n)$ for $n\in\mathbb{Z}$.  
Since $V$ is $C$-stable, 
for sufficiently large $n$ in $\mathbb{N}$ 
we have 
\[
\gamma^n(\alpha(k))=
\alpha(c'_{n}\gamma^n(k){c'}^{-1}_{n})
\in \alpha(c'_{n}V{c'}^{-1}_{n})=\alpha(V)\,.
\]

This shows that 
$\alpha(\mathcal{K}_{\gamma,V})=\mathcal{K}_{\gamma,\alpha(V)}$, 
hence 
$\alpha(\overline{\mathcal{K}}_{\gamma,V})=\overline{\mathcal{K}}_{\gamma,\alpha(V)}$  
for all $C$-stable, compact, open subgroups $V$ of $G$. 

Now, 
if 
$V$ runs through a 
neighborhood base of $e$ 
consisting of $C$-stable compact, open subgroups 
(which exists by Lemma~\ref{lem:cp=>flat}) 
then 
$\alpha(V)$ 
does as well.  
Since 
$K_{\gamma}$ can be defined as 
the intersection of 
the family of all $\overline{\mathcal{K}}_{\gamma,W}$, 
where $W$ runs through a 
%$C$-stable  
neighborhood base of $e$ 
consisting of compact, open subgroups 
as already observed in the proof of Lemma~\ref{lem;tidy_cp-is-aut-stable->cp-stable} 
we obtain $\alpha(K_\gamma)=K_\gamma$.  
\qed 
\vspace{3mm} 

As a corollary of Lemma~\ref{lem:ind-step aut-com-mod-cp} 
we obtain the following result. 

\begin{lemma}\label{lem:almost-abelianG has local tidy subgroups}
Let $\mathcal{H}$ be a group of automorphisms of $G$, 
$C$ a compact subgroup of $G$  
such that 
$[\mathcal{H},\mathcal{H}]$ consists of inner automorphisms in $C$.
Then 
$\mathcal{H}$ 
has local tidy subgroups, 
that is, 
for every finite subset $\mathbf{f}$ of $\mathcal{H}$ 
there is a compact, open subgroup $O$ of $G$ 
such that 
for any $\gamma\in\langle\mathbf{f}\rangle$ 
the group $\gamma(O)$ 
is tidy for each $\alpha\in\mathbf{f}$. 
Moreover, 
$O$ can be chosen $C$-stable. 
\end{lemma}
\pf{}
First, we use induction 
on the cardinality $f\ge 0$ of the finite set $\mathbf{f}$ 
to derive the existence 
of a common tidy subgroup 
for $\mathbf{f}$ 
which is $C$-stable. 
Lemma~\ref{lem:cp=>flat} 
proves the induction hypothesis in the case $f=0$ 
and provides a basis for the induction. 
Assume 
the induction hypothesis 
is already proved for sets of cardinality $f-1\ge 0$, 
and assume that $\mathbf{f}$ is a finite set of cardinality $f$. 
Choose an element $\gamma$ in $\mathbf{f}$ 
and put $\mathcal{A}=\mathbf{f}\setminus\{\gamma\}$. 
Then the induction hypothesis 
implies that there is a common $C$-stable tidy subgroup for $\mathcal{A}$ 
and Lemma~\ref{lem:ind-step aut-com-mod-cp} 
implies that 
the same holds for $\mathbf{f}$.    

If $O$ is a common $C$-stable tidy subgroup for $\mathbf{f}$,  
then 
the first step in the proof of Lemma~\ref{lem:ind-step aut-com-mod-cp} 
shows that $\gamma(O)$ is $\alpha$-tidy 
for each $\alpha$ in $\mathbf{f}$ 
since $\gamma\in\mathcal{H}$. 
\qed
\vspace{3mm}

\pf{Proposition \ref{prop:almost abelian=>flat}}
Finally, an application of Theorem~5.5 in~\cite{tidy<:commAut(tdlcG)} enables us to derive Proposition~\ref{prop:almost abelian=>flat} from Lemma~\ref{lem:almost-abelianG has local tidy subgroups}. 
\qed

% SECTION 4

\section{Topologically simple, non-linear groups of arbitrary flat rank }\label{sec:cases_equality}
We think that the algebraic and the geometric rank of a Coxeter group of finite rank are equal. 
More generally, Corollary~7 in \cite{cocp-proper-CAT0-S} in conjunction with 
\cite[Theorem~C]{local-struc(lengthS_curv<)} leads us to believe that the rank of a proper $\CAT0$-space on which a group $G$ acts properly discontinuously and cocompactly, will turn out to be a quasi-isometry invariant and will be equal to the algebraic rank of $G$. 

\subsection{A family of non-linear groups with equality of ranks}
\label{subsec:NonLinearFamily}
Looking for examples with $\ark(W)=\rk(|W|_0)$ led us to the existence of non-linear, topologically simple groups of arbitrary flat rank. 

\begin{theorem}
\label{thm:E(non-linear-simple-tdlcG(any_rank))}
For every natural number $n \ge 1$ there is a non-linear, topologically simple, compactly generated, 
locally compact, totally disconnected group of flat rank  $n$. 
\end{theorem}

\pf{}
The idea of the proof is 
to show the existence 
of a sequence 
of connected Coxeter diagrams $(D_n)_{n\in\mathbb{N}}$ 
such that for each $n \in \mathbb{N}$, we have: 
\begin{enumerate}
\item 
if $W_n$ denotes the Coxeter group with diagram $D_n$, 
then $\ark(W_n)=\rk(|W_n|_0)=\dim(|W_n|_0)$; 
\item 
there is a Kac-Moody root datum of associated Coxeter diagram $D_n$, and a finite field $k$ such that the corresponding Kac-Moody group $G_n$ over $k$ is center-free and not linear over any field. 
\end{enumerate}
The required 
non-linear, topologically simple, 
compactly generated, 
locally compact, totally disconnected group 
of flat rank  $n$ 
may then be taken to be the completion $\overline{G}_n$. 

\vspace{3mm} 

We now specify a sequence of Coxeter diagrams satisfying the two conditions above. 
Fix a field $k$ of cardinality at least $5$. 
Let $D_1$ be a cycle of length $r=5$ all of whose edges 
are labeled $\infty$. 
Then $W_1$ is Gromov-hyperbolic 
and therefore satisfies 
the first condition above. 
There is 
a Kac-Moody root datum 
whose Coxeter diagram is $D_1$ 
and 
whose group of $k$-points 
satisfies 
the conditions of \cite[Theorem~4.C.1]{top-simpl&com-superrig&non-lin(Kac-MoodyG)} 
(take the remark succeeding that theorem into account). 
This shows that 
the second condition above 
is satisfied as well.  
This settles the case $n=1$. 

If $n>1$ 
let $D_n$ be the diagram 
obtained by 
joining every vertex of $D_1$ 
to every vertex of 
a diagram of type $\widetilde{A}_n$ 
by an edge labeled $\infty$. 
The translation subgroup 
of the special subgroup 
corresponding to 
the $\widetilde{A}_n$-subdiagram   
is an abelian subgroup 
of rank $n$ 
in $W_n$. 
Furthermore,  
the dimension of the $\CAT0$-realization 
of $W_n$ 
can be seen to be $n$ as well 
since the size of 
the maximal spherical subsets of $S$ 
does not grow. 
As the dimension of $|W_n|_0$ 
is an upper bound 
for $\rk(|W_n|_0)$, 
the first condition above 
is satisfied for $D_n$ 
for any $n>1$. 
For every $n>1$ 
the Kac-Moody data 
chosen in the case $n=1$ 
can be extended to a Kac-Moody data 
such that 
the associated Coxeter group 
has diagram $D_n$.  
(This amounts to 
extend 
the combinatorial data 
--- 
that is, 
the finite index set 
of cardinality $n$, 
the generalized $n\times n$ Cartan matrix, 
the choice of $n$ of vectors each in  
a free $\mathbb{Z}$-module  
of rank $n$ 
and its dual,  
such that 
the matrix of the 
pairing 
between these two sets of vectors 
is the given generalized Cartan matrix 
---
defining 
the Kac-Moody group functor. 
The most restrictive of these tasks is 
to arrange 
for the coefficients 
of a general Cartan matrix 
to yield the given Weyl group.   
This is possible 
as soon as 
the edge labels of 
the Coxeter diagram  
of the Weyl group 
belong to the set $\{2, 3, 4, 6, \infty\}$.)
By our choice then %This implies that 
$G_n$ contains $G_1$, 
and the second condition 
is also satisfied 
for any $n>1$. 
The group $G_n$ is center-free whenever the above described Kac-Moody root datum is chosen to be the adjoint datum for the given generalized Cartan matrix. 
\qed\vspace{3mm}

\subsection{Algebraic rank, after Krammer}
\label{subsec:quotingKrammer}
We next observe that thanks to a theorem of Krammer's, the algebraic rank of a Coxeter group of finite rank 
can be computed from its Coxeter diagram. 
In order to state that theorem, we first introduce the notion of standard abelian subgroup of a Coxeter group. 

\begin{definition}
\label{def:standard-ab<CoxeterG}
Let $(W,S)$ be a Coxeter system. 
Let $I_1,\ldots , I_n \subseteq S$ be irreducible, non-spherical and pairwise perpendicular.
For any $i$, let $H_i$ be a subgroup of $W_{I_i}$ defined as follows. 
If $I_i$ is affine, $H_i$ is the translation subgroup; otherwise, $H_i$ is any infinite cyclic subgroup of $W_{I_i}$. 
The group $\prod_i H_i$ is called a \emph{standard free abelian subgroup}.
\end{definition}

The algebraic rank of a standard abelian subgroup $\prod_i H_i$ as above equals 
$\sum_{I_i\ \text{affine}} (\#I_i-1)+\sum_{I_i\ \text{not affine}} 1$. 
Moreover all possible choices of the subsets $I_1,\ldots , I_n \subseteq S$ can be enumerated, so the maximal algebraic rank of a Coxeter group is achieved by some standard free abelian subgroup because of the following theorem \cite[6.8.3]{thesisKrammer}. 
 
\begin{theorem}
\label{thm:alg-rk(CoxeterG)} 
Let $W$ be an arbitrary Coxeter group of finite rank. 
Then any free abelian subgroup of $W$ has a finite index subgroup which is conjugate to a subgroup of some standard free abelian subgroup.
\qed
\end{theorem}

In particular, the algebraic rank of a Coxeter group of finite rank can be computed from its Coxeter diagram. 

\subsection{Application to isomorphism problems}
\label{subsec:isomorphism}
We finish by mentioning an application of the notion of flat rank to isomorphisms of discrete groups, via the use of super-rigidity. 

\begin{proposition}
\label{prop:isomorphism}
Let $\Lambda$ be 
an almost split Kac-Moody group 
over a finite field $\mathbb{F}_q$ of characteristic $p$. 
We assume that $q\geq 4$ 
and $W({1 \over q})<\infty$, 
where $W(t)$ is the growth series of the Weyl group $W$. 
Let ${\bf H}$ be a simple algebraic group 
defined and isotropic over a local field $k$, archimedean or not. 
Assume there exists a surjective group homomorphism from $\Lambda$ to a lattice of ${\bf H}(k)$.
Then, we have: ${\rm char}(k)=p$ and $\frk(\overline\Lambda)=k$-${\rm rk}({\bf H})$. \end{proposition}

\pf{}
Let $\Gamma$ be a lattice in ${\bf H}(k)$. 
Let $\eta\colon \Lambda\to\Gamma$ be a surjective group homomorphism. 
Since ${\bf H}$ is $k$-isotropic, ${\bf H}(k)$ is non compact and $\Gamma$ is infinite. 
By surjectivity of $\eta$, we have $\eta\bigl(Z(\Lambda)\bigr)\leqslant Z(\Gamma)$, 
so by factoring out modulo the centers 
we obtain a surjective group homomorphism 
$\bar\eta\colon \Lambda/Z(\Lambda)\to\Gamma/Z(\Gamma)$. 
By Borel density \cite[II.4.4]{EdM3.17} 
$Z(\Gamma)$ lies in the finite center $Z({\bf H})$, 
so $\Gamma/Z(\Gamma)$ is infinite. 
By the normal subgroup property for $\Lambda$ 
\cite{factor+normal<Thm(lattices<prod(G)),int(ind-cocyles(KMG))}, 
any proper quotient of $\Lambda/Z(\Lambda)$ is finite. 
This implies that $\bar\eta$ is an isomorphism, 
so that the group $\Lambda/Z(\Lambda)$ is linear in characteristic $p$. 
By \cite[Proposition 4.3]{class+non-lin-prop(KM-lattices)}, this implies ${\rm char}(k)=p$. 
In particular, $k$ is a non-archimedean local field, and ${\bf H}(k)$ is totally disconnected. 

Denoting by $\widetilde{\bf H}$ the adjoint quotient of ${\bf H}$ 
and applying the embedding theorem of \cite[3.A]{top-simpl&com-superrig&non-lin(Kac-MoodyG)}, 
we obtain a closed embedding of topological groups 
$\mu\colon \overline\Lambda\hookrightarrow \widetilde{\bf H}(k)$. 
This is where we use that $q\geq 4$ and $W({1 \over q})<\infty$. 
The closed subgroup $\mu(\overline\Lambda)$ of $\widetilde{\bf H}(k)$ contains the lattice $\Gamma$. 
All groups in the chain $\Gamma\leqslant \mu(\overline\Lambda)\leqslant \widetilde{\bf H}(k)$ 
being unimodular, 
$\mu(\overline\Lambda)$ is of finite covolume in $\widetilde{\bf H}(k)$. 
Therefore by \cite[Main Theorem]{strong-approx(ssG(functionF))}, $\mu(\overline\Lambda)$ contains $\widetilde{\bf H}(k)^+$. 
In view of the remark following Corollary \ref{cor:flatrk(semi-simpleG)}, this implies $\frk(\overline\Lambda)=k$-${\rm rk}({\bf H})$. 
\qed
\vspace{3mm}

Replacing the target group ${\bf H}(k)$ by a closed isometry group of a suitable $\CAT0$-space, this result might be generalized to provide restrictions on the existence of discrete actions on $\CAT0$-cell complexes. 
For this, one shall replace the commensurator super-rigidity by super-rigidity of irreducible lattices, as proved in \cite{superrig(irred-latt)+geom-splitting}. 

Assume now that the target lattice $\Gamma$ is replaced by an infinite Kac-Moody lattice (and that both Kac-Moody groups are split). 
Using again the normal subgroup property \cite{factor+normal<Thm(lattices<prod(G)),int(ind-cocyles(KMG))}, we are led to an isomorphism $\bar\eta:\Lambda/Z(\Lambda)\simeq\Gamma/Z(\Gamma)$. 
In this situation, a much more precise answer than ours was given by Caprace and M\"uhlherr 
\cite{iso(KMG)-preserve-bound<G}. 
Their proof is of geometric and combinatorial nature (Bruhat-Tits fixed point lemma and twin root datum arguments). 

%	\bibliographystyle{alpha}
%	\bibliography{%
%	short-Abk,%
%	local,%
%	AAGruppen,%
%	%Abk,%
%	Baeume,% 
%	%GGtheory,% 
%	Gitter,% 
%	NCRaeume,% 
%	%arithmetic,% 
%	buildings,% 
%	diverses,% 
%	%growth,%
%	%measures,% 
%	%qi,% 
%	%quartett,% 
%	%short-Abk,% 
%	%tilings,% 
%	topG%, 
%	%types,% 
%	%unclassified%
%	}
%	\end{document}

\end{document}